\documentclass[11pt]{article} 
\usepackage{amsmath,amssymb,amsfonts,amsthm,amscd,graphicx,psfrag,epsfig}
\usepackage{stmaryrd}
\usepackage[titletoc,toc]{appendix}
\usepackage{hyperref}
\hypersetup{colorlinks, linkcolor=blue}
\usepackage[all]{xy}           
\usepackage{marginnote}  
\usepackage{float}
\usepackage{scalerel}
\usepackage{amsmath,amssymb,amsfonts,amsthm,amscd,graphicx,psfrag,epsfig,mathabx}
\usepackage[titletoc,toc]{appendix}

\newtheorem{theorem}{Theorem}[section]

\newtheorem{theorem-definition}[theorem]{Theorem-Definition}
\newtheorem{theorem-construction}[theorem]{Theorem-Construction}
\newtheorem{lemma-definition}[theorem]{Lemma--Definition}
\newtheorem{lemma-construction}[theorem]{Lemma--Construction}
\newtheorem{lemma}[theorem]{Lemma}
\newtheorem{proposition}[theorem]{Proposition}
\newtheorem{corollary}[theorem]{Corollary}
\newtheorem{conjecture}[theorem]{Conjecture}
\newtheorem{definition}[theorem]{Definition}

\newenvironment{remark}[1][Remark.]{\begin{trivlist}
\item[\hskip \labelsep {\bfseries #1}]}{\end{trivlist}}

\newcommand{\old}[1]{}

\newcommand{\Q}{{\mathbb Q}}

\newcommand{\G}{{\rm G}}

\renewcommand{\P}{{\mathbb P}}

\newcommand{\lra}{\longrightarrow}

\newcommand{\be}{\begin{equation}}
\newcommand{\ee}{\end{equation}}
\newcommand{\bt}{\begin{theorem}}
\newcommand{\et}{\end{theorem}}
\newcommand{\bd}{\begin{definition}}
\newcommand{\ed}{\end{definition}}
\newcommand{\bp}{\begin{proposition}}
\newcommand{\ep}{\end{proposition}}

\newcommand{\bl}{\begin{lemma}}
\newcommand{\el}{\end{lemma}}
\newcommand{\bc}{\begin{corollary}}
\newcommand{\ec}{\end{corollary}}
\newcommand{\bcon}{\begin{conjecture}}
\newcommand{\econ}{\end{conjecture}}
\newcommand{\la}{\label}

\newcommand{\bi}{\begin{itemize}}
\newcommand{\ei}{\end{itemize}}
\newcommand{\bs}{\begin{split}}
\newcommand{\es}{\end{split}}
\begin{document}

\date{}
 
\title{The Strong Suslin Reciprocity Law}

\author{Daniil Rudenko}

\maketitle

\tableofcontents

\section{Introduction}
\subsection{Search for motivic complexes}
For an arbitrary field $F$ and integers $n$ and $m$  Beilinson   (\cite{B82}, \cite{B87}) suggested a definition of motivic cohomology groups
$H^n(F,\mathbb{Q}(m))$ as appropriate pieces of higher algebraic $K-$theory groups, defined earlier by Quillen. In several cases ($n=m$ and $n=1, m=2$) Suslin proved that these groups could be identified with the cohomology groups of very explicit complexes (\cite{Su85}, \cite{Su91}). For the case $n=m$ this implies an isomorphism between motivic cohomology and Milnor $K-$theory.

One may hope that this construction can be generalized to arbitrary integers $n$ and $m.$ This leads to a definition of very explicit complexes (polylogarithmic motivic complexes), whose cohomology groups conjecturally coincide with rational motivic cohomology groups, see \cite{G91a}.
The beauty and elegance of this construction comes with a price: proving  functoriality in $F$ is very hard. Difficulties already exist in the case of Milnor $K-$theory, which we will discuss in detail below. 

The main goal of this paper is to establish existence of norm maps for the cohomology groups of polylogarithmic motivic complexes for the case $n=m-1,$ namely the degree ``next to Milnor K-theory." This implies some corollaries, including  the {\it Strong Suslin Reciprocity Law} conjectured by Goncharov (\cite{G05}, p. 53, Conj. 6.2). We will also formulate an application to hyperbolic scissors congruence theory, which was the author's original motivation behind this project.

\subsection{Main definitions and results}

Let $A$ be an abelian group. Denote by $A_\mathbb{Q}=A\otimes_{\mathbb{Z}} \mathbb{Q}$ its rationalization. 

Let $F$ be an arbitrary field. We will associate  several $\mathbb{Q}-$vector spaces with $F$.  
We start with $F^{\times}_{\mathbb{Q}}:=F^{\times}\otimes_{\mathbb{Z}} \mathbb{Q}:$ the multiplicative group of a field, made rational.  We have an exact sequence
$$
0 \lra  R_1(F) \lra \mathbb{Q}[F]_1 \lra F^\times_{\mathbb{Q}} \lra 0,
$$
where $\mathbb{Q}[F]_1$ is a $\mathbb{Q}-$vector space freely generated by vectors $[p]_1$ for $p\in \P^1_F= F\cup \{\infty\}$  subject to the relations $[0]_1=[1]_1=[\infty]_1=0$.  The subscript $1$ in the symbol $\mathbb{Q}[F]_1$ is only for bookkeeping and refers to the motivic weight. The subspace $R_1(F)$ is generated by elements $[x_1,x_2]_1=[x_1]_1+[x_2]_1-[x_1x_2]_1.$ These symbols satisfy the following {\it cocycle} relation:
\be
\begin{split}
&[x_1,x_2]_1+[x_1 x_2,x_3]_1=[x_1,x_3]_1+[x_1 x_3,x_2]_1.\\
\end{split}
\ee 

Next, we define the {\it second Bloch group}\footnote{The name {\it  Bloch group} was used by Suslin \cite{Su91} for the kernel of the map 
$\delta: B_2(F) \lra \Lambda^2 F^\times$ discussed below. 
We use the terminology of \cite{G91a}, where the group $B_2(F)$ was called the Bloch group, and its higher analogs ${\cal B}_n(F)$ were introduced and called {\it higher Bloch groups}.}  $B_2(F)_{\mathbb{Q}}$ as the cokernel of a map
$$
R_2(F) \lra \mathbb{Q}[F]_2. 
$$
Here the vector space $\mathbb{Q}[F]_2$  is defined just as $\mathbb{Q}[F]_1$  is defined and the space $R_2(F)$ is freely generated by the symbols $[x_1,x_2]_2,$ whose image in $\mathbb{Q}[F]_2$  equals the sum 
$$[x_1]_2+[x_2]_2+[x_3]_2+[x_4]_2+[x_5]_2,$$ where the terms $x_i$ satisfy the following $5-$periodic recurrence:
$$
x_{n+1}=\frac{1-x_n}{x_{n-1}}.
$$
Explicitly,
$$
x_3=\frac{1-x_2}{x_1},\ \ \ \ \ x_4=\frac{x_1+x_2-1}{x_1x_2}, \ \ \ \ \ x_5=\frac{1-x_1}{x_2}.
$$
For some values of $x_1$ and $x_2$ in the formulas above an indeterminate term of the form $0/0$ or $\infty / \infty$ may appear. In this case we omit such terms in the formula for the map $R_2(F)\lra \mathbb{Q}[F]_2.$ We will denote the image of $[x]_2\in \mathbb{Q}[F]_2$ in $B_2(F)_{\mathbb{Q}}$  by $\{x\}_2.$ It is easy to see that the following relations hold in $B_2(F)_{\mathbb{Q}}:$
\be
\begin{split}
& \{a\}_2+\{1-a\}_2=0,\\
&\{a\}_2+\left \{ \frac{1}{a} \right \}_2=0 \\
&\{a\}_2- \{b\}_2+\left \{ \frac{b}{a} \right \}_2-\left \{ \frac{1-a^{-1}}{1-b^{-1}}\right \}_2+\left \{ \frac{1-a}{1-b} \right \}_2=0.\\
\end{split}
\ee 
The third relation is called the {\it five-term relation} and is often used in defining $B_2(F).$

\begin{remark}
The definition above is motivated by the fact that for $F=\mathbb{C}$ the map
$$
\mathcal{L}_2 \colon B_2(\mathbb{C}) \lra \mathbb{R} 
$$ 
sending $[x]_2$ to the {\it Bloch-Wigner dilogarithm} $\mathcal{L}_2(x)$ is well-defined, thanks to Abel's equation for dilogarithm, see \cite{G94}. We omit the rationalization sign, because $B_2(\mathbb{C})$ is uniquely divisible, see \cite{Su91}, Remark 5.1.
\end{remark}

The following complex  ${\cal B}(F,2)$ is called {\it the weight two polylogarithmic complex} :
$$
B_2(F)_{\mathbb{Q}} \stackrel{\delta}{\lra} \Lambda^2 F^\times_{\mathbb{Q}},
$$
where $\delta \{x\}_2=(1-x)\wedge x.$ The fact that $\delta$ annihilates the five-term relation can be checked by a direct computation:
$$
\delta\left (\sum_{i=1}^5 \{x_i\}_2\right )=\sum_{i=1}^5 (1-x_i) \wedge x_i=\sum_{i=1}^5 (x_{i-1}x_{i+1})\wedge x_i=\sum_{i=1}^5 x_{i-1}\wedge x_i - \sum_{i=1}^5 x_{i}\wedge x_{i+1}=0.
$$

Conjectures of  Goncharov (\cite{G91a}, p.5, Conj. B) imply that the vector spaces $F^{\times}_{\mathbb{Q}}$ and $B_2(F)_{\mathbb{Q}}$ are the first two graded components of a conjectural graded Lie coalgebra $\mathcal{L}_{F},$ {\it the Lie coalgebra of mixed Tate motives over $F$.} The category of finite-dimensional graded co-representations of $\mathcal{L}_{F}$ should be equivalent to a conjectural abelian category of mixed Tate motives over $F.$  The cohomology of the Chevalley-Eilenberg complex of $\mathcal{L}_{F}$ in each graded component should coincide with rational motivic cohomology groups of the field $F$:
$$
H^{i}([\Lambda^{\bullet} \mathcal{L}_{F}]_j,\mathbb{Q})=H_{\mathcal{M}}^{i}(F, \mathbb{Q}(j)).
$$ 
This suggests the existence of a very explicit presentation for motivic cohomology groups. The above is known for $i=j=n>0:$ the corresponding cohomology of $\mathcal{L}_{F}$ is equal to the cokernel of the map 
$$
B_2(F)_{\mathbb{Q}} \otimes \Lambda^{n-2} F^{\times}_{\mathbb{Q}} \stackrel{\delta}{\lra} \Lambda^n F^{\times}_{\mathbb{Q}},
$$
which coincides with the rationalized Milnor $K-$theory of $F$. It follows from the results of Suslin \cite{Su85} that 
$$
K_n^M(F) \otimes \mathbb{Q}=H_{\mathcal{M}}^{n}(F, \mathbb{Q}(n)).
$$
Results of Goncharov \cite{G91b} suggest that the third graded component of $\mathcal{L}_{F}$ should be generated by elements $[x]_3 \in \mathbb{Q}[F]_3$ for $x \in \mathbb{P}^1_F$ subject to the relations representing functional equations for the single-valued version of the {\it classical trilogarithm function} $\mathcal{L}_3,$ defined in \cite{G91a}. This motivates the introduction of {\it truncated polylogarithmic complexes} ${\cal B}_2(F,n)$:
$$
 [B_2(F)_{\mathbb{Q}} \otimes \Lambda^{n-2} F^{\times}_{\mathbb{Q}}]_t  \stackrel{\delta}{\lra} \Lambda^n F^{\times}_{\mathbb{Q}},
$$
where the first term is defined as the cokernel of the  map
$$
\mathbb{Q}[F]_3 \otimes \Lambda^{n-3} F^{\times}_{\mathbb{Q}} \stackrel{\delta}{\lra} B_2(F) \otimes \Lambda^{n-2} F^{\times}_{\mathbb{Q}}
$$
sending $[x]_3\otimes y_3 \wedge \ldots \wedge y_n$ to $\{x\}_2\otimes x \wedge y_3 \wedge \ldots \wedge y_n.$ 

\begin{conjecture}[(Goncharov)]
For an arbitrary field $F$ and $n\geq 2$ the kernel of the map 
$$
\delta \colon [B_2(F)_{\mathbb{Q}} \otimes \Lambda^{n-2} F^{\times}_{\mathbb{Q}}]_t \lra \Lambda^n F^{\times}_{\mathbb{Q}}
$$
is isomorphic to the motivic cohomology group $H^{n-1}_{\mathcal{M}}(F, \mathbb{Q}(n)).$
\end{conjecture}
This is known for $n=2$ thanks to the work of Suslin \cite{Su91}. We will denote this kernel by 
$H^{n-1}_{\G}(F, \mathbb{Q}(n)).$
Our main result is:
\begin{theorem}\label{main}
Let $k$ be an arbitrary field and $n\geq 3$. Then the following sequence is exact:
$$
0 \lra H^{n-1}_{\G}(k, \mathbb{Q}(n)) \lra H^{n-1}_{\G}(k(t), \mathbb{Q}(n)) \stackrel{\oplus \partial_P}{\lra} \bigoplus \limits_{P\in \mathbb{A}^1_k} H^{n-2}_{\G}(k_P, \mathbb{Q}(n-1)) \lra 0.
$$	
\end{theorem}
Here $k_P$ is the  residue field at the point $P \in \mathbb{A}^1_k.$ The {\it residue maps} $\partial_P$ will be defined below. We will formulate and prove a similar statement for  Milnor $K-$theory, which was used by Bass, Tate \cite{BT} and Suslin \cite{Su79}  to construct norm maps in Milnor $K-$theory. Their construction works in our case as well:
\begin{corollary}\label{C1}
Let $F$ be a perfect field. For every finite field extension $L/F$ there is a norm homomorphism
$$
Nm_{L/F} \colon H^{n-1}_{\G}(L, \mathbb{Q}(n)) \lra H^{n-1}_{\G}(F, \mathbb{Q}(n))
$$
that makes the functor $H^{n-1}_{\G}(-, \mathbb{Q}(n))$ into a Rost cycle module as defined in \cite{R}.
\end{corollary}
This implies the following conjecture of Goncharov (\cite{G05}, p. 53, Conjecture 6.2):
\begin{corollary}[(Strong Reciprocity Law)]\label{C2}
Let $X$ be a compact smooth algebraic curve over $\mathbb{C}.$ For every  $n\geq 3$  the map
$$
Res \colon  {\cal B}_2(\mathbb{C}(X),n)\lra {\cal B}_2(\mathbb{C},n-1)
$$
is null-homotopic. Here  the total residue map $Res$ is a sum over all points $P\in X$ of local residue maps $\partial_P.$
\end{corollary}
\begin{remark}
For curves of genus $g\leq 1$ Goncharov  found a formula for a contracting homotopy  (\cite{G05}, Theorems 6.5, 6.14). For curves of higher genus explicit construction of $h$ remains unknown.
\end{remark}
To formulate the next corollary, we recall some notions from the theory of scissors congruences. Let $\mathcal{P}(\mathbb{H}^3)$ be a $\mathbb{Q}-$vector space generated by classes of hyperbolic polytopes modulo the scissors congruence equivalence relation. The so-called Dehn invariant 
$$
D\colon \mathcal{P}(\mathbb{H}^3) \lra  \mathbb{R} \otimes \mathbb{R}/ 2 \pi \mathbb{Z}
$$
associates to a polytope  the sum $\sum_{e\in E} l_e \otimes \alpha_e,$ where $E$ is the set of edges of a polytope, $l_e$ is the length of the corresponding edge and $\alpha_e$ is the corresponding dihedral angle. There is a chain map $S[i]$
\be 
\begin{gathered}
    \xymatrix{
    &  B_2(\mathbb{C})\ar[d]^{S[1]} \ar[r] & \Lambda^2\mathbb{C}^{\times}_{\mathbb{Q}}\ar[d]^{S[2]}\\
    &  \mathcal{P}(\mathbb{H}^3) \ar[r] &  \mathbb{R} \otimes \mathbb{R}/ 2 \pi \mathbb{Z}
    }
 \end{gathered},
 \ee 
where map $S[1]\colon B_2(\mathbb{C}) \lra \mathcal{P}(\mathbb{H}^3)$
sends an element $\{z\}_2$ to the scissors congruence class of the ideal tetrahedron with vertices $\infty,0,1, z$ and 
\be
\begin{split}
&S[2]\bigl((1-z) \wedge z\bigr)=-(1-z) \wedge z+(1-\bar{z}) \wedge \bar{z}=\\ 
&2\left(|z|\otimes \arg(1-z)-|1-z|\otimes \arg(z) \right)\in \mathbb{R} \otimes \mathbb{R}/ 2 \pi \mathbb{Z}.
\end{split}
\ee

\begin{corollary}\label{Scis}
Let $X$ be a  smooth projective curve over $\mathbb{C}.$ Then there is a map
$$
 h\colon \Lambda^3 \mathbb{C}(X)^{\times} \lra \mathcal{P}(\mathbb{H}^3)
$$
such that for any three meromorphic functions $f_1, f_2, f_3$ the Dehn invariant of a polytope $h(f_1 \wedge f_2 \wedge f_3)$ is equal to  $(S[2 \circ Res)(f_1 \wedge f_2 \wedge f_3)$ and the hyperbolic volume of the polytope $h(f_1 \wedge f_2 \wedge f_3)$ is given by the following convergent integral:
$$
\frac{1}{2\pi i}\int_{X(\mathbb{C})} r_2(f_1,f_2,f_3), 
$$
where the form $r_2(f_1,f_2,f_3)$ equals to 
$$
\frac{1}{6}\sum_{\sigma \in \mathbb{S}_3} (-1)^{\sigma}\Bigl ( \log|f_{\sigma(1)}| d \log|f_{\sigma(2)}| \wedge d  \log|f_{\sigma(3)}|- 3 \log|f_{\sigma(1)}|  d \arg(f_{\sigma(2)}) \wedge d \arg(f_{\sigma(3)})\Bigr ).
$$
\end{corollary}

\begin{remark}
	The statement of the Corollary \ref{Scis} was the original motivation behind this project. It was prompted by  scissors congruence properties of {\it Schl\"afli orthoschemes}, discovered by J.-P. Sydler 	\cite{Sy}. Hyperbolic orthoschemes appear in the image of the map $h$ from Corollary \ref{Scis} for $X=\mathbb{P}^1.$  

\end{remark}

\subsection{The structure of the article}
Towards proving Theorem \ref{main}, we start by discussing the classical proof of an exact sequence for Milnor $K-$theory: Theorem \ref{BT}. There is nothing original in our exposition, but the proof of Theorem \ref{main} is based on it. This proof is based on the properties of a certain filtration, which we will use extensively in the proof of Theorem \ref{main}. We do not write in detail proofs of Corollaries \ref{C1}, \ref{C2}  because they are not different from the Milnor $K-$theory case, which is explained in \cite{Su79}. Corollary \ref{Scis} follows from results of \cite{G05} and \cite{Du}.

Our proof of the main result occupies Sections \ref{S3}-\ref{S9}. In Section \ref{S3} one can find the proof modulo two key statements: Lemmas  \ref{L1} and \ref{L2}. Proofs of \ref{L1} and \ref{L2} are based on the study of the co-residue map, which are defined in the derived category only.  The most non-trivial step in the proof of Theorem \ref{main} is to show that these maps are well-defined, which is done in  Sections $\ref{S6}$ and $\ref{S7}.$ 

\subsection{Notations and conventions}\label{SecConv}

All our results are valid only modulo torsion, so we work everywhere with $\mathbb{Q}-$vector spaces. Thus it is convenient to omit the rationalization sign from the notation, which we do everywhere, starting with Section \ref{Sec2}, except Section \ref{SecSus}. In particular, we will use $F^{\times}$ instead of $F^{\times}_{\mathbb{Q}},$ $K_n^M(F)$ instead of $K_n^M(F)_{\mathbb{Q}},$ $B_2(F)$ instead of $B_2(F)_{\mathbb{Q}},$ and so on.

We use notation $\mathbb{S}^k V$ and $\Lambda^k V$ for symmetric and wedge powers of a vector space. Tensor product $V\otimes W$ is understood as tensor product over $\mathbb{Q}.$

Next, for  $n\geq 3$ we denote the cokernel of the map 
$$
\mathbb{Q}[F]_3 \otimes \Lambda^{n-3} F^{\times} \stackrel{\delta}{\lra} B_2(F) \otimes \Lambda^{n-2} F^{\times}
$$
 by the symbol $[B_2(F)\otimes \Lambda^{n-2} F^{\times}]_t$. For every sub-factor space $X$ of the space $B_2(F)\otimes  \Lambda^{n-2} F^{\times}$ we will denote by the symbol $[X]_t$ its projection to $[B_2(F)\otimes \Lambda^{n-2} F^{\times}]_t$.

A few remarks about filtered vector spaces are in order. By a {\it filtered vector space $(V,\mathcal{F})$}   we mean a vector space $V$ equipped with   an increasing filtration 
$$
\mathcal{F}_0 V \subset \mathcal{F}_1 V \subset \ldots \subset  V.
$$
We will denote by $gr_k^{\mathcal{F}}V$ the associated graded space $\mathcal{F}_k V / \mathcal{F}_{k-1} V.$ A morphism 
$$
f \colon (V,\mathcal{F}) \lra (W,\mathcal{G})
$$ 
of filtered vector spaces is called {\it strictly compatible with filtration} if 
$$
f(\mathcal{F}_k V )=f(V) \cap \mathcal{G}_k W
$$
Every subspace or factor space of a filtered vector space inherits the filtration. Similarly, there exists a natural filtration on the tensor product of a finite number of filtered vector spaces. A detailed exposition of these notions can be found in \cite{D}.

\subsection{Acknowledgement}
I would like to thank S. Gorchinskiy, A. Levin, O. Martin and A. Goncharov  for their invaluable help with preparing the manuscript. Also I thank Gerhard Paseman for his help on a draft of this article.

\section{Results of Bass, Tate and Milnor} \label{Sec2}

\subsection{Milnor K-theory}

We will start with discussing a well-known  analogue of Theorem \ref{main} for Milnor $K-$theory. It was proved by Bass and Tate  \cite{BT} for $n=2$  and by Milnor \cite{M} in general.

Let $F$ be an arbitrary field. For every integer $n$ we define $K_n^M(F)$ -- the {\it Milnor K-group  of $F$} to be the quotient of the abelian group
$\Lambda^n F^{\times}$ by the subgroup generated by tensors $(1-x)\wedge x \wedge y_3 \wedge \ldots \wedge y_n.$  The image of an element $y_1 \wedge y_2 \wedge \ldots \wedge y_n$ in $K_n^M(F)$ is denoted by $\{y_1, y_2,\ldots,y_n\}.$ Our definition agrees with the classical one modulo $2-$torsion. In accordance with the convention from Section \ref{SecConv}, we denote by $K_n^M(F)$ the $\mathbb{Q}-$vector space $K_n^M(F)_{\mathbb{Q}}.$

Our next goal is to define residue maps. Suppose that $\nu$ is a discrete valuation of $F$ with residue field $F_{\nu}$ and uniformizer $f_\nu.$    First we define a map
\be
\partial_\nu \colon \Lambda^n F^{\times} \lra  \Lambda^{n-1} F_{\nu}^{\times}.
\ee
For any elements $u_1, \ldots, u_n \in F^\times$ with vanishing discrete valuation we let
$$
\partial_\nu(u_1\wedge u_2 \wedge \ldots \wedge u_n)=0
$$
and 
$$
\partial_\nu(f_\nu\wedge u_2 \wedge \ldots \wedge u_n)=\overline{u_2}\wedge \ldots \wedge \overline{u_n}.
$$
Here $\bar{u}$ is the residue of $u$ in $F_\nu^{\times}.$ This map can be uniquely extended by linearity to $\Lambda^n F^{\times}.$ It is easy to see that it induces a well-defined map from $K_n^M(F)$ to $K_{n-1}^M (F_{\nu}).$

We will be especially interested in the case when $F=k(t)$ and  the valuation $\nu_P$ corresponds to a point $P\in \mathbb{A}^1_k$. We will denote the residue field  by  $k_P,$ and a uniformizer by $f_P.$  The uniformizer will always be chosen to be a monic polynomial.

\begin{theorem}[(Milnor, Tate)] \label{BT}
 The following sequence is exact for $n\geq 2$:
$$
0 \lra K_n^M(k) \stackrel{j}{\lra} K_n^M(k(t))\stackrel{\oplus \partial_P}{\lra} \bigoplus \limits_{P\in \mathbb{A}^1_k} K_{n-1}^M (k_P) \lra 0.
$$
\end{theorem}
Here the map $K_n^M(k) \stackrel{j}{\lra}  K_n^M(k(t))$
is induced by an inclusion of $k$ in $k(t).$ 

\subsection{Some corollaries of Theorem \ref{BT}}

Assume that $F$ is perfect. Theorem \ref{BT} can be used to construct a norm map 
$$
Nm_{L/K}\colon K_n^M(L) \lra K_n^M(F)
$$
for any finite field extension $L/F.$ For this consider any point $P\in \mathbb{A}^1_F$ with residue field $L$ (such a point exists by the primitive element theorem). Thanks to Theorem \ref{BT}, for any element $w\in K_n^M(L)$ there exists an element $W \in K_{n+1}^M(F(t))$ with residue (at the point $P$) equal to $w$  and with vanishing residues at all other points. Then one can define $Nm_{L/F}(w)$ to be the residue of $W$ at $\infty \in \mathbb{P}^1.$ 

One can show that this construction does not depend on the choice of $P.$ The proof is straightforward for the case when one works modulo torsion, as we do. Using the norm map, one can establish all usual properties of norms. In particular, we obtain the following nontrivial result:
\begin{corollary}[(Suslin Reciprocity Law)]
	Let $X$ be a compact smooth algebraic curve over $\mathbb{C}.$ Then for $n\geq 1$ the map
$$
Res \colon  K_n^M(\mathbb{C}(X))\lra K_{n-1}^M(\mathbb{C})
$$
is  equal to zero. Here  the total residue map $Res$ is a sum over all points $P\in X$ of local residue maps $\partial_P.$
\end{corollary}

 All the details of these constructions can be found in \cite{Su79}.

\subsection{Proof of  Theorem \ref{BT}}  \label{Sec2.3}
\begin{proof}

Denote the $\mathbb{Q}-$vector space $k(t)^{\times}$ by $D.$ It is naturally graded by the degree d:
$D_0=k^{\times}$ and for $d$ greater than zero  $D_d$ is freely generated by monic irreducible polynomials of degree $d.$ We will use interchangeably the notation $D_{\leq d}$ and $\mathcal{F}_d D$ for $\bigoplus \limits_{i \leq d} D_i.$ The filtration $\mathcal{F}_\bullet$ can be extended to all wedge powers of $k(t)^{\times}$ and thus defines a filtration on $ K_n^M(k(t)).$ Explicitly, $\mathcal{F}_d K_n^M(k(t))$ is generated by symbols $\{f_1,\ldots, f_n\},$ where $f_i$ are irreducible polynomials of degree $\leq d.$ 

Let  $P \in \mathbb{A}^1_k$ be a point of degree $d.$ Then the map $\partial_P$ vanishes on all symbols  $\{g_1,\ldots, g_n\},$ where $g_i$ are irreducible polynomials of degree  $<d.$ We deduce that the following map is well defined:
$$
gr_d^{\mathcal{F}} K_n^M(k(t))\stackrel{\oplus \partial_P}{\lra} \bigoplus \limits_{\deg(P)=d} K_{n-1}^M (k_P).
$$ 
We claim that this map is an isomorphism. To see that this map $\oplus \partial_P$ is surjective, consider any element $\{r_1, r_2, \ldots, r_{n-1}\} \in K_{n-1}^M (k_P).$ For each $i$ let $\widetilde{r_i}$ be a polynomial of minimal degree with residue $r_i.$ Then
$$
\partial_P \{f_P, \widetilde{r_1}, \widetilde{r_2}, \ldots, \widetilde{r}_{n-1}\} =\{r_1, r_2, \ldots, r_{n-1}\}. 
$$  
If $Q$ is any other point of degree $d,$ 
$$
\partial_Q \{f_P, \widetilde{r_1}, \widetilde{r_2}, \ldots, \widetilde{r}_{n-1}\}=0.
$$
This implies that the map $\oplus \partial_P$ is surjective.

To see that $\oplus \partial_P$ is injective, for every point $P$ we construct {\it co-residue} maps 
$$
c_P \colon K_{n-1}^M (k_P) \lra gr_d^{\mathcal{F}} K_n^M(k(t)).
$$ 
The definition is straightforward: for every element $\{r_1, r_2, \ldots, r_{n-1}\} \in K_{n-1}^M (k_P)$ we put 
 $$
 c_P(\{r_1, r_2, \ldots, r_{n-1}\})=\{f_P, \widetilde{r_1}, \widetilde{r_2}, \ldots, \widetilde{r}_{n-1}\} \in gr_d^{\mathcal{F}} K_n^M(k(t)).
 $$
\begin{lemma}
\label{Lemma2.3}
The co-residue map 
$$
c_P \colon K_{n-1}^M (k_P) \lra gr_d^{\mathcal{F}} K_n^M(k(t))
$$
is well-defined.
 \end{lemma}
 
 \begin{proof}
First, we show that the map $c_P$ is linear in each variable. Consider a pair of  elements in $K_{n-1}^M (k_P):$
$$
S_1=\{s_1, r_2, \ldots, r_{n-1}\}, \ \ 
S_2=\{s_2, r_2, \ldots, r_{n-1}\}.
$$ 
Then,
$$
S_1+S_2=\{s_1s_2, r_2, \ldots, r_{n-1}\}.
$$
It is easy to see that, 
 \be
 \begin{split}
 & c_P(S_1)=\{f_P, \widetilde{s_1}, \widetilde{r_2}, \ldots, \widetilde{r}_{n-1}\}, \\
 & c_P(S_2)=\{f_P, \widetilde{s_2}, \widetilde{r_2}, \ldots, \widetilde{r}_{n-1}\}, \\
 & c_P(S_1+S_2)=\{f_P, \widetilde{s_1s_2}, \widetilde{r_2}, \ldots, \widetilde{r}_{n-1}\},\\
 \end{split}
 \ee
 so
 \be
 c_P(S_1)+c_P(S_2)-c_P(S_1+S_2)=\left \{ f_P, \dfrac{\widetilde{s_1}\widetilde{s_2}}{\widetilde{s_1s_2}},\widetilde{r_2}, \ldots, \widetilde{r}_{n-1} \right \}.
 \ee
 Note that from the division algorithm it follows that there exists a polynomial $q(t)$ such that
 $$
 \widetilde{s_1}\widetilde{s_2}=q  f_P+\widetilde{s_1s_2},
 $$
 so that
  $$
 \frac{\widetilde{s_1}\widetilde{s_2}}{\widetilde{s_1s_2}}+\frac{-q  f_P}{\widetilde{s_1s_2}}=1.
 $$
This implies
 \be
 \begin{split}
 & c_P(S_1)+c_P(S_2)-c_P(S_1+S_2)= \\
 & \left \{ f_P, \dfrac{\widetilde{s_1} \widetilde{s_2}}{\widetilde{s_1s_2}},\widetilde{r_2}, \ldots, \widetilde{r}_{n-1} \right \}= \\
 &\left \{ f_P, 1-\frac{-q f_P}{\widetilde{s_1s_2}},\widetilde{r_2}, \ldots, \widetilde{r}_{n-1} \right \}=\\
 &\left \{ \frac{-q f_P}{\widetilde{s_1s_2}}, 1-\frac{-q f_P}{\widetilde{s_1s_2}},\widetilde{r_2}, \ldots, \widetilde{r}_{n-1} \right \}-\left \{ \frac{-q}{\widetilde{s_1s_2}},\dfrac{\widetilde{s_1} \widetilde{s_2}}{\widetilde{s_1s_2}} ,\widetilde{r_2}, \ldots, \widetilde{r}_{n-1} \right \}=\\
  &-\left \{ \frac{-q}{\widetilde{s_1s_2}},\dfrac{\widetilde{s_1} \widetilde{s_2}}{\widetilde{s_1s_2}} ,\widetilde{r_2}, \ldots, \widetilde{r}_{n-1} \right \}.\\
 \end{split}
 \ee
  Since 
 $$
 \deg(q)\leq \max\left(\deg(\widetilde{s_1})+\deg(\widetilde{s_2})-\deg(f_P), 0 ) \right)<d,
 $$
 the expression 
$
-\left \{ \dfrac{-q}{\widetilde{s_1s_2}},\dfrac{\widetilde{s_1} \widetilde{s_2}}{\widetilde{s_1s_2}} ,\widetilde{r_2}, \ldots, \widetilde{r}_{n-1} \right \}
$
 lies in $\mathcal{F}_{d-1} K_n^M(k(t)),$ so it vanishes in $gr_d^{\mathcal{F}} K_n^M(k(t)).$
Thus $c_P(S_1)+c_P(S_2)=c_P(S_1+S_2).$ 

It is easy to see that $c_P$ is antisymmetric and vanishes on elements 
$\{r_1, 1-r_1, r_3,\ldots, r_{n-1}\},$ so this map is well-defined.
\end{proof}

We claim that the map $\sum c_P$ is the inverse of $\oplus \partial_P.$
Obviously $\partial_P \circ c_P=id$ and  $\partial_P \circ c_Q=0$ if points $P$ and $Q$ are distinct. It remains to show that 
\be \label{EqCor}
\sum c_P \partial_P=id.
\ee
Notice that this equality holds for elements in $gr_d^{\mathcal{F}}K_n^M \left(k(t)\right)$ of the form 
$\{f_P,\widetilde{r_1}, \widetilde{r_2}, \ldots, \widetilde{r_{n-1}}\}$
where $\widetilde{r_i}$ are elements of degree less than $d.$

\begin{lemma}
Vector space $gr_d^{\mathcal{F}} K_n^M (k(t))$ is generated by elements 
$$\{f,g_2, g_3, \ldots, g_{n}\},$$
for polynomials $f$ and $g_i$ such that $\deg(f)=d$ and $\deg(g_i)<d.$
\end{lemma}
\begin{proof}
Denote by $\mathcal{G}_s gr_d^{\mathcal{F}} K_n^M (k(t))$ a subspace of $gr_d^{\mathcal{F}} K_n^M (k(t))$ generated by elements
$$
\{f_1,f_2,\ldots, f_k, g_{k+1}, \ldots, g_{n}\}\in gr_d^{\mathcal{F}} K_n^M (k(t)),
$$
where $k\leq s,$ polynomials $f_i$ are monic irreducible of degree $d$ and the degree of polynomials $g_i$ is less than $d.$ Clearly, 
$$
gr_d^{\mathcal{F}} K_n^M (k(t))=\mathcal{G}_n gr_d^{\mathcal{F}} K_n^M (k(t)).
$$
To prove the lemma it is sufficient to show that 
$$
\mathcal{G}_s gr_d^{\mathcal{F}} K_n^M (k(t))=\mathcal{G}_{s-1} gr_d^{\mathcal{F}} K_n^M (k(t))
$$
for $2\leq s \leq n.$

We  present a proof for $n=2,$ the general case being similar.
Consider two distinct monic polynomials $f_1, f_2$ of degree $d.$ Then 
$$0=\left \{ 1-\frac{f_1}{f_2} ,\frac{f_1}{f_2} \right \}=\{f_1,f_2\}-\{f_1,f_2-f_1\}+\{f_2,f_2-f_1\}.$$
Note that since $f_2$ and $f_1$ are monic, $\deg(f_2-f_1)<d,$ so
$$
\{f_1,f_2\} \in \mathcal{G}_1 gr_d^{\mathcal{F}} K_2^M (k(t))
$$ From this the statement follows.
\end{proof}
We conclude that the map
$$
gr_d^{\mathcal{F}} K_n^M(k(t))\stackrel{\oplus \partial_P}{\lra} \bigoplus \limits_{\deg(P)=d} K_{n-1}^M (k_P)
$$ 
is an isomorphism, from which Theorem \ref{BT} easily follows.
\end{proof}

\section{The plan of the proof of Theorem \ref{main}}\label{S3}
To give a precise statement of Theorem  \ref{main}, we need to extend the definition of the residue map $\partial_P$ to truncated polylogarithmic complexes.
\subsection{Residue map}
Let $F$ be a field with discrete valuation $\nu$ and residue field $F_{\nu}.$ Denote by $f_\nu$ a uniformizer. Our goal is to define a chain map $\partial_{\nu}$
\be 
\begin{gathered}
    \xymatrix{
   &  [B_2(F) \otimes \Lambda^{n-2}F^{\times}]_t \ar[d]^{\partial_{\nu}} \ar[r] & \Lambda^{n}F^{\times} \ar[d]^{\partial_{\nu}} \\
   &  [B_2(F_{\nu}) \otimes \Lambda^{n-3}F_{\nu}^{\times}]_t  \ar[r]& \Lambda^{n-1}F_{\nu}^{\times}  }
                          \end{gathered}
 \ee 
for $n\geq 3.$
We have already defined this map on  $\Lambda^{n}F^{\times}$ in the previous section. 
On 
$$
\{x\}_2\otimes  y_3 \wedge \ldots \wedge  y_{n} \in [B_2(F) \otimes \Lambda^{n-2}F^{\times}]_t
$$ 
we define it by the formula 
$$
\{\overline{x}\}_2 \otimes \partial_{\nu}(y
_3\wedge \ldots \wedge  y_{n}) \in  [B_2(F_{\nu}) \otimes \Lambda^{n-3}F_{\nu}^{\times}]_t.
$$
It is easy to see that $\partial_{\nu}$ is a chain map of complexes
$$
{\cal B}_2(F,n)\lra {\cal B}_2(F_{\nu},n-1).
$$
Theorem \ref{main} is equivalent to the statement that 
the complex of chain maps 
$$
0\lra {\cal B}_2(k,n)\lra {\cal B}_2(k(t),n)\lra \bigoplus_{P \in \mathbb{A}^1_k}{\cal B}_2(k_P,n-1)\lra 0
$$
is exact on cohomology.

\subsection{Filtration by degree.}
The main tool in the proof of Theorem \ref{BT} was an auxiliary filtration $\mathcal{F}_{\bullet}$ on $K_n^M(k(t))$ induced by the degree filtration on $k(t)^{\times}$.
This filtration can be extended to the complex ${\cal B}_2(k(t),n).$ For this we first define a  filtration $\mathcal{F}_{\bullet}$ on the Bloch group $B_2(k(t))$ as a pre-image of the filtration  $\mathcal{F}_{\bullet}$ on $\Lambda^2 k(t)^{\times}$ 
under $\delta.$ In the next section we will use results of Suslin to describe this filtration more explicitly. This defines a filtration on $B_2(k(t)) \otimes \Lambda^{n-2}k(t)^{\times},$ which descends to $\left [ B_2(k(t)) \otimes \Lambda^{n-2}k(t)^{\times}\right ]_t.$

Recall that a map 
$$
f \colon (V,\mathcal{F})\lra (W,\mathcal{F})
$$
of filtered spaces is called {\it strictly compatible with the filtration} if 
$$
f(\mathcal{F}_i V)=f(V) \cap \mathcal{F}_i W.
$$
\begin{lemma}
The map
	$$
	\delta \colon [B_2(k(t)) \otimes \Lambda^{n-2}k(t)^{\times}]_t \lra \Lambda^{n}k(t)^{\times}, 
	$$
	is strictly compatible with the filtration $\mathcal{F}.$
\end{lemma}
\begin{proof}
	The statement follows from the proof of Theorem \ref{BT}.
\end{proof}

\begin{lemma}
To prove Theorem \ref{main} it is sufficient to show that for every positive degree $d$  the map
\be \label{eqlem3}
\oplus \partial_P \colon gr^{\mathcal{F}}_d{\cal B}_2(k(t),n)\lra \bigoplus_{\deg(P)=d}{\cal B}_2(k_P,n-1)
\ee
is a quasi-isomorphism.
\end{lemma}
\begin{proof}
Assume that (\ref{eqlem3}) is a quasi-isomorphism.
We will use the fact that $\mathcal{F}_0B_2(k(t))=B_2(k),$ which will be proved in the next section.
Our goal is to show that the complex
\be \label{seqlem3}
0 \lra H^{n-1}_{\G}(k, \mathbb{Q}(n)) \lra H^{n-1}_{\G}(k(t), \mathbb{Q}(n)) \stackrel{\oplus \partial_P}{\lra} \bigoplus \limits_{P\in \mathbb{A}^1_k} H^{n-2}_{\G}(k_P, \mathbb{Q}(n-1)) \lra 0
\ee
is exact for $n\geq 3$, where $H^{n-1}_{\G}(F, \mathbb{Q}(n))$ is the kernel of the map  
$$
	\delta \colon [B_2(F) \otimes \Lambda^{n-2}F^{\times}]_t \lra \Lambda^{n}F^{\times}.
$$
First, consider an element $x \in H^{n-1}_{\G}(k(t), \mathbb{Q}(n))$ in the kernel of $\oplus\partial_P.$ Then 
$$
x\in \mathcal{F}_0 [B_2(k(t)) \otimes \Lambda^{n-2}k(t)^{\times}]_t,
$$
by (\ref{eqlem3}). But 
	$$
\mathcal{F}_0[B_2(k(t)) \otimes \Lambda^{n-2}k(t)^{\times}]_t=
[\mathcal{F}_0(B_2(k(t))) \otimes \Lambda^{n-2}\mathcal{F}_0 (k(t)^{\times})]_t=[B_2(k) \otimes \Lambda^{n-2}k^{\times}]_t,
$$
so $x \in H^{n-1}_{\G}(k, \mathbb{Q}(n)).$ So the complex (\ref{seqlem3}) is exact in the middle term.
To show that $\oplus \partial_P$ is surjective, notice that $\delta$ is strictly compatible with $\mathcal{F},$ so
$gr_d^{\mathcal{F}} H^{n-1}_{\G}(k(t), \mathbb{Q}(n))$
is the kernel of the map 
$$
	\delta \colon gr_d^{\mathcal{F}} [ B_2(k(t)) \otimes \Lambda^{n-2}k(t)^{\times} ]_t \lra gr_d^{\mathcal{F}} \Lambda^{n}k(t)^{\times}.
$$
It follows from the quasi-isomorphism (\ref{eqlem3}) that this kernel is isomorphic to 
$$
\bigoplus \limits_{P\in \mathbb{A}^1_k} H^{n-2}_{\G}(k_P, \mathbb{Q}(n-1)).
$$
From this the lemma follows.
\end{proof}

\subsection{Filtration by support}
It remains to prove that the  map
$$
\oplus \partial_P \colon gr^{\mathcal{F}}_d{\cal B}_2(k(t),n) \lra
\bigoplus_{\deg(P)=d}{\cal B}_2(k_P,n-1)
$$
is a quasi-isomorphism. For this we introduce another filtration on 
$
gr^{\mathcal{F}}_d{\cal B}_2(k(t),n)
$:
{\it the filtration by support}. This is an increasing filtration $\mathcal{G}_1\subset \mathcal{G}_2\subset \ldots\subset \mathcal{G}_n.$
On $gr_d^{\mathcal{F}} \left(\Lambda^n k(t)^{\times}\right)$ we define the filtration by placing tensors 
$$
f_1  \wedge \ldots \wedge f_s \wedge g_{s+1} \wedge \ldots \wedge g_{n}
$$
with polynomials $f_i, g_i$ such that $\deg(f_i)=d,
$ and $\deg(g_i)<d$ in $\mathcal{G}_s.$  It  is easy to see that 
$$
gr_s^\mathcal{G}gr_d^{\mathcal{F}} \left(\Lambda^n k(t)^{\times}\right)=\Lambda^s D_d \otimes \Lambda^{n-s} D_{<d}.
$$
On $gr^{\mathcal{F}}_d(B_2(k(x)))$ we define the filtration as the pre-image of the filtration $\mathcal{G}$ on $gr_d^{\mathcal{F}} \left(\Lambda^2 k(t)^{\times}\right)$ under $\delta.$ It will be computed explicitly in the next section. Finally, on $gr^{\mathcal{F}}_d  [B_2(k(t)) \otimes \Lambda^{n-2}k(t)^{\times} ]_t$ the filtration is obtained by projecting the tensor product of the corresponding filtrations on the components.

\subsection{What remains to be proved}

Theorem \ref{main} follows from Lemmas \ref{L1} and \ref{L2}.

\begin{lemma}
	The complex 
	$$
	gr_1^\mathcal{G}gr_d^{\mathcal{F}}{\cal B}_2(k(t),n)
	$$
	is quasi-isomorphic to 
	$$ 
\bigoplus_{\deg(P)=d}{\cal B}_2(k_P,n-1).
$$
\label{L1}
\end{lemma}
Lemma \ref{L1} will be proved in Section \ref{S8} using a construction of  {\it co-residue} maps $c_P$. The difficulty is that the co-residue maps can be defined  in the derived category only. 
\begin{lemma}
	For $s>1$ the complex 
	$$
	gr_s^\mathcal{G}gr_d^{\mathcal{F}}{\cal B}_2(k(t),n)
	$$
	is acyclic.
	\label{L2}
\end{lemma}

 Lemma \ref{L2} is not very hard and will be proved in Section \ref{S9}.
To see that Lemmas \ref{L1} and \ref{L2} imply Theorem \ref{main} we also need to show that the map 
	$$
	\delta \colon gr_d^{\mathcal{F}} \left [ B_2(k(t)) \otimes \Lambda^{n-2}k(t)^{\times} \right ]_t \lra gr_d^{\mathcal{F}} \Lambda^{n}k(t)^{\times}
	$$
is strictly compatible with filtration $\mathcal{G}$; compatibility will follow from the proof of Lemma \ref{L2}.

\section{Anatomy of the group $B_2(k(t))$}\label{S4}
\subsection{Introduction}
The goal of this section is to describe a set of generators for the group $B_2(k(t))$. The results of this section can be formulated very explicitly (see Corollary \ref{Cor}) but the proofs use the connection with higher algebraic $K-$theory, established by Suslin. 

Recall from Section \ref{Sec2.3} that for $d\geq 1$ the vector space $D_d\subset k(t)^{\times}$ is freely generated by monic irreducible polynomials of degree $d.$ Also, $D_0=k^{\times}.$ 

\begin{definition}
For a point $P\in \mathbb{A}^1_k$ of degree $d$ consider a linear map
$$
\rho_P \colon D_{<d}\lra k_P^{\times},
$$
which sends a polynomial $f\in k(t)^{\times}$ of degree less than $\deg(P)$ to its residue $\bar{f}.$ We denote the kernel of the map $\rho_P$ by $B_P.$
\end{definition}

The map $\rho_P$ is surjective: any residue class $r \in k_{P}^{\times}$ is the image of a unique polynomial of degree less than $\deg(P),$ which we denote by $\widetilde{r}.$ 

Consider a map 
$$
\alpha_P \colon R_1(k_{P}) \lra B_P,
$$ 
defined on generators $[r_1,r_2]_1$ by the formula 
$$
\alpha_P ([r_1,r_2]_1)= \dfrac{\widetilde{r_1} \widetilde{r_2}}{\widetilde{r_1 r_2}}
$$ 
and extended to $R_1(k_{P})$ by linearity. Clearly, this map is well-defined and surjective.

\begin{remark}
We will later use  the fact that $D_d$ is rather close to the rational group ring of $k_P^{\times}$ and $B_P$ is related to the square of the augmentation ideal of the group ring. More precisely, there exists a surjective map
$$
\mathbb{Q}[k_P^{\times}]\lra D_{< d},
$$ 	
sending $r$ to $\tilde{r}.$ The kernel is generated by elements $[r_1]+[r_2]-[r_1r_2],$ where $\deg(\widetilde{r_1})+\deg(\widetilde{r_2})<d.$
\end{remark}

\subsection{Result  of  Suslin}\label{SecSus}
The second cohomology group of the weight two polylogarithmic complex ${\cal B}^\bullet(F,2)$ is naturally identified with a rationalized Milnor $K-$group. Suslin proved\footnote{Suslin's results are more precise, since he does not neglect torsion. To emphasize the difference, in this section we do not omit the rationalization sign, as we do in the rest of the text. Also, Suslin's results are proved under an assumption that the field is infinite. For a finite fields all the results of this section hold automatically, since all the groups involved  vanish. }  in \cite{Su91} that its first cohomology group is naturally identified with the indecomposable part of $K_3(F)_{\mathbb{Q}}$. Hence, the following sequence is exact:
$$
0 \lra K_3^{ind}(F)_\Q \lra  B_2(F)_{\mathbb{Q}} \stackrel{}{\lra}   \Lambda^2 F^\times_{\mathbb{Q}} \lra K_2^{M}(F)_\Q \lra 0.
$$
We will use two other results from $K-$theory coming from the relation between the $K-$theory of $k$ and $k(t).$ First, the following exact sequence is a special case of Theorem \ref{BT}:
$$
0 \lra K_2^{M}(k) \lra K_2^{M}(k(t))   \lra \bigoplus \limits_{P\in \mathbb{A}^1_k} k_P^{\times} \lra 0.
$$
Second, from the localization sequence and $\mathbb{A}^1-$homotopy invariance  for algebraic $K-$theory it follows that the  embedding
$$
K_3^{ind}(k) \lra K_3^{ind}(k(t))
$$
is an isomorphism.
Combining these two results with Suslin's theorem, we get the following statement.
\begin{corollary} \label{Sus}
	The following sequence is exact:
$$
0 \lra B_2(k)_{\mathbb{Q}}\lra  B_2(k(t))_{\mathbb{Q}} \lra   \dfrac{\Lambda^2 k(t)^\times_{\mathbb{Q}}}{\Lambda^2 k^\times_{\mathbb{Q}}} \lra \bigoplus \limits_{P\in \mathbb{A}^1_k} \left( k_P^{\times} \right)_{\mathbb{Q}} \lra 0.
$$
\end{corollary}

\subsection{Degree filtration}

The associated graded factor $gr_d^\mathcal{F}(\Lambda^2 k(t)^\times)$ can be described explicitly: 
$$
gr_d^\mathcal{F}(\Lambda^2 k(t)^\times)= \left ( D_d \otimes D_{<d} \right ) \oplus \Lambda^2 D_d.
$$
The vector space $B_2(k(t))$ carries a filtration induced from the filtration on $\Lambda^2 k(t)^\times,$ which we will also denote by $\mathcal{F}.$ It follows from Corollary \ref{Sus} that $\mathcal{F}_0 B_2(k(t))=B_2(k)$ and that the following sequence is exact for $d>0$:
 $$
0 \lra gr_d^\mathcal{F}[ B_2(k(t))] \lra   gr_d^\mathcal{F}[\Lambda^2 k(t)^\times] \lra \bigoplus \limits_{\deg(P)=d} k_P^{\times} \lra 0.
$$

\begin{lemma}
The following sequence is exact:	
$$
0 \lra \bigoplus \limits_{\deg(P)=d} B_P \lra gr_d^\mathcal{F}[ B_2(k(t))]  \lra \Lambda^2 D_d \lra 0.
$$
\end{lemma}

\begin{proof} The projection $gr_d^{\mathcal{F}}[ B_2(k(t))]  \lra \Lambda^2 D_d$ is surjective. Indeed, given two irreducible monic polynomials $f$ and $g$ of degree $d,$ consider the symbol
$\left \{ \dfrac{f}{g}\right \}_2 \in B_2(k(t)).$
Since 
$$
\delta \left \{ \dfrac{f}{g}\right \}_2=\left ( 1-\dfrac{f}{g} \right )\wedge\dfrac{f}{g} =f\wedge g +g \wedge (g-f)+(g-f)\wedge f,
$$
and $g-f$ has degree less than $d,$ it follows that  $\left \{ \dfrac{f}{g}\right \}_2 $ projects to  $f \wedge g$ in $\Lambda^2 D_d.$ 
Consider the following commutative diagram with exact rows and surjective columns:
\be 
\begin{gathered}
    \xymatrix{
    & 0  \ar[r]& gr_d^{\mathcal{F}}[ B_2(k(t))] \ar[d] \ar[r] & gr_d^{\mathcal{F}}[\Lambda^2 k(t)^\times]  \ar[d]^{} \ar[r]^{}& \bigoplus \limits_{\deg(P)=d} k_P^{\times} \ar[d]^{}\ar[r] & 0\\
    & 0 \ar[r]& \Lambda^2 D_d  \ar[r]& \Lambda^2 D_d  \ar[r] &0 \ar[r] &0}
                          \end{gathered}.
 \ee 
Applying the snake lemma, we get an exact sequence
$$
0 \lra Ker \left ( gr_d^{\mathcal{F}}[ B_2(k(t))]  \lra \Lambda^2 D_d \right)\lra D_d \otimes D_{<d}  \lra \bigoplus \limits_{\deg(P)=d} k_P^{\times} \lra 0.
$$
On the other hand, $B_P$ was defined as the kernel of the map $D_{<d} \stackrel{\rho_P}{\lra} k_P^{\times},$ so 
$$ 
Ker \left ( gr_d^{\mathcal{F}}[ B_2(k(t))] \lra \Lambda^2 D_d \right)
$$ 
is isomorphic to $\bigoplus \limits_{\deg(P)=d} B_P.$

\end{proof}

Let $P\in \mathbb{A}^1_k$ be a point of degree $d.$ A composition of the map 
$$
\alpha_P \colon R_1(k_{P}) \lra B_P
$$  	
and an embedding 
$$
B_P \hookrightarrow gr_d^\mathcal{F}[ B_2(k(t))]
$$
is a map 
$$
\beta_P\colon R_1(k_{P}) \lra gr_d^\mathcal{F}[ B_2(k(t))].
$$
More explicitly, $\beta_P([r_1,r_2]_1)=\left \{ \dfrac{\widetilde{r_1}\widetilde{r_2}}{\widetilde{r_1r_2}} \right\}_2.$

\begin{remark}
The map $\beta_P$ is well-defined by construction. This is equivalent to the following relation in $gr_d^\mathcal{F}[B_2(k(t))]:$ 
$$
\beta_P([r_1,r_2]_1)+\beta_P([r_1 r_2,r_3]_1)=\beta_P([r_1,r_3]_1)+\beta_P([r_1r_3,r_2]_1).
$$
This relation could also be derived independently by applying the $5-$term relation several times.
\end{remark}

\begin{corollary}\label{Cor}
The vector space $B_2(k(t))$ is generated over its subspace $B_2(k)$ by the following two types of elements:

1) symbols
$
\left \{ \dfrac{f(t)}{g(t)}\right \}_2,
$
 where $f$ and $g$ are irreducible monic polynomials of the same degree;

2) symbols
$
\beta_P([r_1,r_2]_1)=\left \{ \dfrac{\widetilde{r_1}\widetilde{r_2}}{\widetilde{r_1r_2}} \right\}_2,
$
where $r_1$ and $r_2$ are two elements in $k_P^{\times}$ for some point $P.$
\end{corollary}

\begin{remark}
Though this corollary is formulated in a completely elementary way, the author was not able to derive it from the $5-$term relation without using algebraic $K-$theory. If $k$ is algebraically closed this is known as the Rogers dilogarithm identity, see Theorem 8.14 in \cite{Du}. 
\end{remark}

\section{Co-residue map in the derived category}\label{S5}
Our goal is to prove Lemma \ref{L1}, claiming that the complex 
	$$
	gr_1^\mathcal{G}gr_d^{\mathcal{F}}{\cal B}_2(k(t),n)
	$$
	is quasi-isomorphic to 
	$$ 
\bigoplus_{\deg(P)=d}{\cal B}_2(k_P,n-1).
$$
For this we will construct  a co-residue chain map from a complex, which is quasi-isomorphic to ${\cal B}_2(k_P,n-1),$ to $gr_d^{\mathcal{F}}{\cal B}_2(k(t),n)$ with image in $\mathcal{G}_1$. The main difficulty will be to prove the fact that this map is well-defined.

\subsection{Free resolution of the truncated polylogarithmic complex}
To construct a free resolution of the polylogarithmic complex
$$
[B_2(F)\otimes \Lambda^{n-2} F^{\times}]_t \lra \Lambda^n F^{\times}, 
$$
recall the presentation of $B_2,$ $F^{\times}$ and $B_2 \otimes_a F^{\times}$ by generators and relations:
\be
\begin{split}
&0 \lra  R_1(F) \lra \mathbb{Q}[F]_1 \lra F^\times \lra 0,\\
&R_2(F) \lra \mathbb{Q}[F]_2 \lra B_2(F) \lra 0,\\
&\mathbb{Q}[F]_3 \otimes \Lambda^{n-3}F^\times \lra B_2(F)\otimes \Lambda^{n-2} F^\times \lra  [B_2(F)\otimes \Lambda^{n-2} F^\times]_t \lra 0.\\
\end{split} 
\ee
Using standard properties of wedge powers, one obtains the following resolution of $\Lambda^n F^{\times}:$
$$
\mathbb{S}^2 R_1(F) \otimes \Lambda^{n-2} \mathbb{Q}[F]_1 \lra R_1(F) \otimes \Lambda^{n-1} \mathbb{Q}[F]_1 \lra \Lambda^n \mathbb{Q}[F]_1 \lra \Lambda^n F^\times \lra 0
$$
and the following resolution of $[B_2(F)\otimes \Lambda^{n-2} F^\times]_t:$
\be
\begin{split}
&\mathbb{Q}[F]_3 \otimes \Lambda^{n-3}\mathbb{Q}[F]_1\\
&\ \ \ \ \ \ \ \  \ \ \ \ \oplus \\
&\mathbb{Q}[F]_2 \otimes R_1(F) \otimes \Lambda^{n-3} \mathbb{Q}[F]_1\lra \mathbb{Q}[F]_2 \otimes \Lambda^{n-2} \mathbb{Q}[F]_1 \lra   [B_2(F)\otimes \Lambda^{n-2} F^\times ]_t \lra 0.\\
&\ \ \ \ \ \ \ \ \ \ \  \ \oplus\\
&R_2(F) \otimes \Lambda^{n-2} \mathbb{Q}[F]_1\\
\end{split}
\ee
The map $\delta=\delta[0]\colon [B_2(F)\otimes \Lambda^{n-2} F^\times ]_t\lra \Lambda^{n} F^\times$
can be lifted to a map $\delta[i]$ between the resolutions above. Clearly, we have to define
$$
\delta[1] \colon \mathbb{Q}[F]_2 \otimes \Lambda^{n-2} \mathbb{Q}[F]_1 \lra \Lambda^{n} \mathbb{Q}[F]_1
$$
by the formula 
$$
\delta[1] \Big ( [x]_2\otimes [y_3]_1 \wedge \ldots \wedge [y_n]_1 \Big)=[1-x]_1 \wedge [x]_1 \wedge[y_3]_1 \wedge \ldots \wedge [y_n]_1.
$$
The map $\delta[2]$ can be defined in many different ways. 
We put $\delta[2]$ equal to zero on $\mathbb{Q}[F]_3 \otimes \Lambda^{n-3}\mathbb{Q}[F]_1.$
Next, define $\delta[2]$ on $\mathbb{Q}[F]_2 \otimes R_1(F) \otimes \Lambda^{n-3} \mathbb{Q}[F]_1$
by the formula
$$
\delta[2] \Big ( [x]_2\otimes [y_1,y_2]_1 \otimes [z_3]_1 \wedge \ldots \wedge [z_n]_1 \Big)=[y_1,y_2]_1 \otimes [1-x]_1 \wedge [x]_1 \wedge[z_3]_1 \wedge \ldots \wedge [z_n]_1.
$$
Finally, we need to define $\delta[2]$ on $R_2(F) \otimes \Lambda^{n-2} \mathbb{Q}[F]_1.$ It is sufficient to define $\delta[2]$ on $R_2(F).$

For that recall the most simple proof of the five-term relation. The element $[x_1,x_2]_2\in R_2(F)$ is mapped to the element $\sum_{i=1}^5 [x_i]_2\in \mathbb{Q}[F]_2,$
where $x_i\in F$  satisfy the following $5-$periodic recurrence:
$$x_{i+1}=\frac{1-x_i}{x_{i-1}}.$$
Then 
$$\delta \left (\sum_{i=1}^5[x_{i}]_2\right )=\sum_{i=1}^5(1-x_{i})\wedge x_{i} \in \Lambda^2 F^{\times} .$$
Since $(1-x_{i})\wedge x_{i}=x_{i-1}\wedge x_{i}-x_{i}\wedge x_{i+1},$
the sum vanishes telescopically. 
This suggests the following definition of the map 
$
R_2(F) \stackrel{\delta[2]}{\lra} R_1(F)\otimes \mathbb{Q}[F]_1:
$
$$
\delta[2][x_1,x_2]_2=-\sum_{i=1}^5 [x_{i-1},x_{i+1}]_1 \otimes [x_i]_1.
$$

Finally, denote by  $\widetilde{\cal B}_2(F,n)$ the cone of the map $\delta[\bullet].$
For the convenience of the reader, we list here the terms of the complex $\widetilde{\cal B}_2(F,n)$ in each degree:
\be
\begin{split}
&{\bf Degree \ 1:}\Lambda^n \mathbb{Q}[F]_1 \\
&{\bf Degree \ 2:}\Bigl ( R_1(F)\otimes \Lambda^{n-1} \mathbb{Q}[F]_1 \Bigr )\oplus\Bigl ( \mathbb{Q}[F]_2 \otimes \Lambda^{n-2} \mathbb{Q}[F]_1 \Bigr )\\
&{\bf Degree \ 3:}\Bigl (
\mathbb{S}^2 R_1(F) \otimes \Lambda^{n-2} \mathbb{Q}[F]_1 \Bigr ) \oplus \Bigl ( \mathbb{Q}[F]_3 \otimes \Lambda^{n-3}\mathbb{Q}[F]_1\Bigr )\oplus\\
&\Bigl ( \mathbb{Q}[F]_2 \otimes R_1(F) \otimes \Lambda^{n-3} \mathbb{Q}[F]_1 \Bigr) \oplus \Bigl ( R_2(F) \otimes \Lambda^{n-2} \mathbb{Q}[F]_1 \Bigr ).\\
\end{split}	
\ee
Denote by 
$$
\pi \colon \widetilde{\cal B}_2(F,n) \lra {\cal B}_2(F,n)
$$
the cone chain morphism. We have proven the following statement.
\begin{lemma}\label{L5}
The map $\pi$ induces an isomorphism on the cohomology groups in degrees $1$ and $2.$	
\end{lemma}

\subsection{Co-residue map}
Let $P$ be a point of the degree $d$ in $\mathbb{A}^1_k.$ Our goal is to define a co-residue map $c_P$ in the derived category: 
$$
c_P \colon {\cal B}_2(k_P, n) \lra gr_d^{\mathcal{F}}{\cal B}_2(k(t),n+1)
$$ 
with the property that the composition $\partial_P \circ c_P$ is a quasi-isomorphism. Because of  Lemma \ref{L5}, it is sufficient to  construct a map from   $\widetilde{\cal B}_2(k_P,n)$   to $gr_d^{\mathcal{F}}{\cal B}_2(k(t),n+1),$ which we will also denote by $c_P.$

The definition of $c_P$ is completely straightforward, but checking that it is a chain morphism is nontrivial.
First, define 
$$
c_P[1]\colon \Lambda^n \mathbb{Q}[k_P]_1 \lra gr_d^{\mathcal{F}} \Lambda^{n+1} k(t)^{\times}  
$$
by the formula
$$
[r_1]_1 \wedge [r_2]_1 \wedge \ldots  \wedge [r_n]_1 \lra f_P \wedge \widetilde{r_1} \wedge \widetilde{r_2}\wedge \ldots \wedge \widetilde{r_n}.
$$
Next, define 
$$
c_P[2]\colon R_1(k_P)\otimes \Lambda^{n-1} \mathbb{Q}[k_P]_1 \lra gr_d^{\mathcal{F}}{\cal B}_2(k(t),n+1)$$
by the formula
$$
[r_1,r_2]_1 \otimes [r_3]_1 \wedge \ldots  \wedge [r_{n+1}]_1 \lra  \beta_P([r_1,r_2]_1)\otimes \widetilde{r_3} \wedge  \ldots \wedge \widetilde{r_{n+1}}.
$$
Here $\beta_P([r_1,r_2]_1)=\left \{ \dfrac{\widetilde{r_1} \widetilde{r_2}}{\widetilde{r_1 r_2}} \right \}_2 \in B_P \subset gr_d^{\mathcal{F}} B_2(k(t)).$ 
Finally, define 
$$
c_P[2]\colon R_2(k_P)\otimes \Lambda^{n-2} \mathbb{Q}[k_P]_1 \lra gr_d^{\mathcal{F}}{\cal B}_2(k(t),n+1)
$$
by the formula
$$
[r_1]_2 \otimes [r_2]_1 \wedge \ldots  \wedge [r_n]_1 \lra \{\widetilde{r_1}\}_2 \otimes  f_P \wedge \widetilde{r_2} \wedge  \ldots \wedge \widetilde{r_n}.
$$

To check that $c_P$ is a chain map it is enough to show that the composition  $c_P[2]\circ \delta$ vanishes. Its domain $\widetilde{{\cal B}}_2(k_P,n)[3]$ has four direct summands. On 
$\mathbb{Q}[k_P]_3 \otimes \Lambda^{n-3}\mathbb{Q}[k_P]_1$ the map vanishes by definition of $\delta[2].$ We continue in the next sections.

\subsection{Vanishing of $c_P[2]\circ \delta$ on $\mathbb{Q}[k_P]_2 \otimes R_1(k_P) \otimes \Lambda^{n-3} \mathbb{Q}[k_P]_1$.}
Obviously, it is sufficient to consider the case $n=3.$

 \begin{lemma}\label{L3}
 The following identity holds in $[B_2(F)\otimes \Lambda^2 F^{\times}]_t:$
 $$\{a\}_2 \otimes (1-b) \wedge b-\{b\}_2 \otimes (1-a) \wedge a=0.$$	
 \end{lemma}
 \begin{proof}
Consider the five-term relation
 $$
 \{a\}_2- \{b\}_2+\left \{ \frac{b}{a} \right \}_2-\left \{ \frac{1-a^{-1}}{1-b^{-1}} \right \}_2+\left \{ \frac{1-a}{1-b} \right \}_2=0
 $$	
 and multiply it by $\dfrac{b}{a} \wedge \left ( \dfrac{1-a}{1-b} \right).$ We get the expression above modulo terms $\{x\}_2\otimes x \wedge y=\delta(\{x\}_3\otimes y).$ 
 \end{proof}
For $a=\widetilde{r_1}$ and $b=\dfrac{\widetilde{r_2}\widetilde{r_3}}{\widetilde{r_2r_3}}$ we use Lemma \ref{L3} and get that in $B_2(k(t))\otimes \Lambda^2 k(t)^{\times}$
$$
\{\widetilde{r_1} \}_2 \otimes \left ( 1-\dfrac{\widetilde{r_2}\widetilde{r_3}}{\widetilde{r_2r_3}} \right ) \wedge \dfrac{\widetilde{r_2}\widetilde{r_3}}{\widetilde{r_2r_3}}=\left \{ \dfrac{\widetilde{r_2}\widetilde{r_3}}{\widetilde{r_2r_3}}\right \}_2 \otimes (1-\widetilde{r_1}) \wedge \widetilde{r_1}.
$$
Moreover, using arguments  as in the proof of Lemma \ref{Lemma2.3}, the following equality holds in $gr_d^{\mathcal{F}} \left ( B_2(k(t))\otimes \Lambda^2 k(t)^{\times} \right ):$ 
$$
\{\widetilde{r_1} \}_2 \otimes \left ( 1-\dfrac{\widetilde{r_2}\widetilde{r_3}}{\widetilde{r_2r_3}} \right )\wedge \dfrac{\widetilde{r_2}\widetilde{r_3}}{\widetilde{r_2r_3}}=
\{\widetilde{r_1} \}_2 \otimes f_P \wedge \dfrac{\widetilde{r_2}\widetilde{r_3}}{\widetilde{r_2r_3}},
$$
so
$$\left \{ \dfrac{\widetilde{r_2}\widetilde{r_3}}{\widetilde{r_2r_3}}\right \}_2 \otimes (1-\widetilde{r_1}) \wedge \widetilde{r_1} =
\{\widetilde{r_1} \}_2 \otimes f_P \wedge \widetilde{r_2}+\{\widetilde{r_1} \}_2 \otimes f_P \wedge \widetilde{r_3}-\{\widetilde{r_1} \}_2 \otimes f_P \wedge \widetilde{r_2r_3}.
$$
Equivalently,
$$
(c_P[2]\circ \delta) \left ([r_1]_2\otimes [r_2, r_3]_1\right )=0.
$$

\subsection{Vanishing of $c_P[2]\circ \delta$ on $\mathbb{S}^2 R_1(k_P) \otimes \Lambda^{n-2} \mathbb{Q}[k_P]_1$}\label{S6}

 It is enough to check the vanishing of $c_P[2]\circ \delta$ on $\mathbb{S}^2 R_1(k_P) \otimes \Lambda^{n-2} \mathbb{Q}[k_P]_1$ for $n=2.$  For this we use Lemma \ref{ml1}.
  \begin{lemma}\label{ml1}
 	Let $r_1, r_2, r_3 , r_4$ be some elements in $k_P.$ Then 
 $$
 \left \{ \dfrac{\widetilde{r_1}\widetilde{r_2}}{\widetilde{r_1r_2}}\right \}_2 \otimes \dfrac{\widetilde{r_3}\widetilde{r_4}}{\widetilde{r_3r_4}}+
 \left \{ \dfrac{\widetilde{r_3}\widetilde{r_4}}{\widetilde{r_3r_4}}\right \}_2 \otimes \dfrac{\widetilde{r_1}\widetilde{r_2}}{\widetilde{r_1r_2}}=0
 $$
 in $gr^{\mathcal{F}}_d \left [ B_2(k(t))\otimes k(t)^{\times} \right ]_t.$
 \end{lemma}
 
 This is, probably, the most strange part of the proof. We will prove \ref{ml1} after Lemma \ref{id}.      
\begin{lemma}\label{id}
Let $A$ be an abelian group, and $I<\mathbb{Q}[A]$ be the augmentation ideal. Then the elements of the form $([a_1]+[a_2]-[a_3]-[a_4])^2$ with $a_1 a_2=a_3 a_4 \in A$ generate $\mathbb{S}^2 I^2.$	
\end{lemma}
\begin{proof}
Obviously, elements of the form 
\be
\begin{split}
&([a]+[b]-[ab]-[1])\cdot([c]+[d]-[cd]-[1])\\		
\end{split}
\ee
generate $\mathbb{S}^2 I^2.$	 The statement of the lemma follows from the following identity, which can be checked by a direct computation:
\be
\begin{split}
&-4([a]+[b]-[ab]-[1])\cdot([c]+[d]-[cd]-[1])+\\
&([ab]+[cd]-[bd]-[ac])^2+([ab]+[cd]-[ad]-[bc])^2-([bd]+[ac]-[ad]-[bc])^2-\\
&([a]+[bd]-[d]-[ab])^2-([ad]+[b]-[ab]-[d])^2+ ([a]+[bd]-[ad]-[b])^2-\\ 
&([a]+[bc]-[c]-[ab])^2 - ([ac]+[b]-[ab]-[c])^2 + ([a]+[bc]-[ac]-[b])^2-\\
&([c]+[ad]-[a]-[cd])^2-([ac]+[d]-[cd]-[a])^2+([c]+[ad]-[ac]-[d]) ^2-\\
 &([c]+[bd]-[b]-[cd])^2- ([bc]+[d]-[cd]-[b])^2+ ([c]+[bd]-[bc]-[d])^2+\\
 &2([a]+[b]-[ab]-[1])^2 + 2 ([c]+[d] -[cd]-[1])^2=0.\\
\end{split}
\ee
\end{proof}

\begin{proof}[Proof (of Lemma \ref{ml1})]
We will apply Lemma \ref{id} to the group $A=k_P^{\times}.$  There is a map $\varphi_1$ from the group ring $\mathbb{Q}[k_P^{\times}]$ to $D_{<d},$ sending $[r]$ to $\tilde{r}.$ The image of $I^2$ under this map  is contained in $B_P,$ because the sequence 
$$
0 \lra B_P \lra D_{<d} \stackrel{\rho_P}{\lra} k_P^{\times} \lra 0
$$ 
is exact. Furthermore, 
$$
\varphi_1\left(([r_1]-[1])([r_2]-[1])\right)=\alpha_P([r_1,r_2]_1)=\frac{\widetilde{r_1}\widetilde{r_2}}{\widetilde{r_1r_2}}.
$$ 
Denote the composition of the map $\varphi_1\colon I^2 \lra B_P$ with the embedding $B_P \hookrightarrow gr_d^{\mathcal{F}} B_2(t)$ by $\varphi_2:$
$$
\varphi_2(([r_1]-[1])([r_2]-[1]))=\left \{ \dfrac{\widetilde{r_1}\widetilde{r_2}}{\widetilde{r_1r_2}}\right \}_2.
$$
This map can be used to construct a map 
$$
(\varphi_1\cdot \varphi_2) \colon \mathbb{S}^2I^2 \lra  [B_P \otimes   D_{<d}]_t 
$$
sending $\lambda_1\cdot \lambda_2 \in \mathbb{S}^2I^2$ to $\varphi_2(\lambda_1) \otimes  \varphi_1(\lambda_2)+ \varphi_2(\lambda_2) \otimes \varphi_1(\lambda_1).$

To conclude the proof of Lemma \ref{ml1}, we need to show that the map $(\varphi_1\cdot \varphi_2)$ vanishes on $\mathbb{S}^2I^2$, because the expression 
 $$
 \left \{ \dfrac{\widetilde{r_1}\widetilde{r_2}}{\widetilde{r_1r_2}}\right \}_2 \otimes \dfrac{\widetilde{r_3}\widetilde{r_4}}{\widetilde{r_3r_4}}+
 \left \{ \dfrac{\widetilde{r_3}\widetilde{r_4}}{\widetilde{r_3r_4}}\right \}_2 \otimes \dfrac{\widetilde{r_1}\widetilde{r_2}}{\widetilde{r_1r_2}}
 $$
 equals  to  $ (\varphi_1\cdot \varphi_2)\Bigl(([r_1]-[1])([r_2]-[1]) \cdot ([r_3]-[1])([r_4]-[1]) \Bigr).$ 
 
 It follows from Lemma \ref{id} that it is sufficient to check the vanishing of the map $(\varphi_1\cdot \varphi_2)$ on elements 
 $$
 ([r_1]+[r_2]-[r_3]-[r_4])^2 \in \mathbb{S}^2 I^2,
 $$
where $r_1r_2=r_3r_4$ in $k_P^{\times}.$ Notice that in this case
\be
\begin{split}
&\varphi_1([r_1]+[r_2]-[r_3]-[r_4])=\frac{\widetilde{r_1}\widetilde{r_2}}{\widetilde{r_3}\widetilde{r_4}},\\
&\varphi_2([r_1]+[r_2]-[r_3]-[r_4])=\left \{ \dfrac{\widetilde{r_1}\widetilde{r_2}}{\widetilde{r_1r_2}}\right \}_2-\left \{ \dfrac{\widetilde{r_3}\widetilde{r_4}}{\widetilde{r_3r_4 }} \right \}_2.
\end{split}
\ee
The last expression equals to $\left \{ \dfrac{\widetilde{r_1}\widetilde{r_2}}{\widetilde{r_3} \widetilde{r_4}}\right \}_2.$ To see that it is enough to check that the coproduct 
$$
\delta\left( \left \{ \dfrac{\widetilde{r_1}\widetilde{r_2}}{\widetilde{r_1r_2}}\right \}_2-\left \{ \dfrac{\widetilde{r_3}\widetilde{r_4}}{\widetilde{r_3r_4 }} \right \}_2-\left \{ \dfrac{\widetilde{r_1}\widetilde{r_2}}{\widetilde{r_3} \widetilde{r_4}}\right \}_2\right)
$$
vanishes in $gr_d^{\mathcal{F}}[\Lambda^2 k(t)^{\times}],$ which is obvious.

From this we get that
$$
(\varphi_1\cdot \varphi_2)(([r_1]+[r_2]-[r_3]-[r_4])^2)=2\left \{ \frac{\widetilde{r_1}\widetilde{r_2}}{\widetilde{r_3} \widetilde{r_4}}\right \}_2 \otimes  \frac{\widetilde{r_1}\widetilde{r_2}}{\widetilde{r_3} \widetilde{r_4}}= 0,
$$
which finishes the proof of Lemma \ref{ml1}.
 \end{proof}

\subsection{Vanishing of $c_P[2]\circ \delta$ on $R_2(F) \otimes \Lambda^{n-2} \mathbb{Q}[F]_1$}\label{S7}
Clearly, it is enough to consider the case $n=2.$ Our goal is to show the following:
\begin{lemma}\label{ml2}
For any $r_1, r_2 \in k_P^{\times}$ define $r_{i+1}=\dfrac{1-r_i}{r_{i-1}}, \ i=3,4,5.$ 
Then the expression
\be
\begin{split}
	&(\{\widetilde{r_1}\}_2+\{\widetilde{r_2}\}_2+\{\widetilde{r_3}\}_2+\{\widetilde{r_4}\}_2+\{\widetilde{r_5}\}_2)\otimes f_P+ \\
	&+\left \{ \frac{\widetilde{r_1} \widetilde{r_3}}{1-\widetilde{r_2}} \right \}_2\otimes \widetilde{r_2}+\left \{ \frac{\widetilde{r_2} \widetilde{r_4}}{1-\widetilde{r_3}} \right \}_2\otimes \widetilde{r_3}+\left \{ \frac{\widetilde{r_3} \widetilde{r_5}}{1-\widetilde{r_4}} \right \}_2\otimes \widetilde{r_4}+\left \{ \frac{\widetilde{r_4} \widetilde{r_1}}{1-\widetilde{r_5}} \right \}_2\otimes \widetilde{r_5}+\left \{ \frac{\widetilde{r_5} \widetilde{r_2}}{1-\widetilde{r_1}} \right \}_2\otimes \widetilde{r_1} \\
\end{split}
\ee
vanishes in $gr_d^{\mathcal{F}} \left [B_2(k(t))\otimes k(t)^{\times} \right ]_t.$ 
\end{lemma}
This will finish showing how $c_P[2]\circ \delta$ vanishes, because 
$$
\beta_P([r_i, r_{i+2}]_1)=\left \{ \frac{\widetilde{r_i} \widetilde{r_{i+2}}}{\widetilde{r_i r_{i+2}}}\right \}_2=\left \{ \frac{\widetilde{r_i} \widetilde{r_{i+2}}}{1-\widetilde{r_{i+1}}}\right \}_2.
$$
We start with a lemma whose proof is inspired by Gauss lemma, used in the classical proof of quadratic reciprocity.

\begin{lemma}
\label{Gauss}
For $1\leq i\leq 5$ there exist polynomials $g_i(t), h_i(t)$ of degrees less than $d,$ such that $r_i=\dfrac{\overline{h_i}}{\overline{g_i}}$ in $k_P^{\times}$ and the sequence $\lambda_i=\dfrac{g_i}{h_i}$ satisfies the recurrence
	$$
	\lambda_{i+1}=\frac{1-\lambda_i}{\lambda_{i-1}}.
	$$
\end{lemma}
\begin{proof}
	We will suppose that $d$ is odd: the even case is similar. For the residue $r_1\in k_P$ there exist polynomials $g_1(t), h_1(t)\in k(t)$ of degrees less or equal to $\left \lfloor\frac{d}{2} \right \rfloor,$ such that the residue of $\dfrac{g_1}{h_1}$ equals to $r_1.$ Indeed, consider $k_P$ as a vector space over $k$ with basis $1, \overline{t}, \ldots, \overline{t^{d-1}}$ and denote by $V$ its subspace, spanned by $1, \overline{t}, \ldots, \overline{t^{\lfloor d/2 \rfloor}}.$  Multiplication by $r_1$ is a linear  automorphism of $k_P$ and $\dim(V)>\frac{\dim(k_P)}{2},$ so there exists a vector $g_1 \in V$ such that $r_1  g_1=f_1 \in V.$ 
	 
	 Similarly, we can find $g_2(t), h_2(t)\in k(t)$ of degrees less or equal to $\left \lfloor\frac{d}{2} \right \rfloor,$ such that the residue of $\dfrac{g_2}{h_2}$ equals to $r_2.$ There exists a unique sequence $\lambda_i,$ such that $\lambda_1=\dfrac{g_1}{h_1}$ and $\lambda_2=\dfrac{g_2}{h_2},$ which satisfies the required recurrence. Then 
	 $$
	 \lambda_3=\dfrac{1-\dfrac{g_2}{h_2}}{\dfrac{g_1}{h_1}}=\frac{(h_2-g_2)h_1}{h_2g_1}
	 $$
	 and we define $g_3=(h_2-g_2)h_1$ and $h_3=h_2g_1.$ Similarly,
	 $$
	 \lambda_4=\dfrac{1-\dfrac{g_3}{h_3}}{\dfrac{g_2}{h_2}}=\dfrac{g_1h_2-h_1h_2+g_2h_1}{g_1g_2}=\dfrac{g_4}{h_4}
	 $$
	 and
	$$
	 \lambda_5=\dfrac{1-\dfrac{g_4}{h_4}}{\dfrac{g_3}{h_3}}=\dfrac{(h_1-g_1)h_2}{h_1g_2}=\dfrac{g_5}{h_5}.
	 $$
	 Obviously, $\deg(g_i), \deg(h_i)$ are less than $d$ for $1\leq i\leq 5$ as required. 
\end{proof}

To simplify the notation a little bit we will omit \ $\tilde{}$ \ and \ $\bar{}$\ \ symbols till the end of Section \ref{S7}. 
First we want to rewrite $\bigl(\sum_{i=1}^{5}\{r_i\}_2 \bigr)\otimes f_P $ as an element of $[B_P \otimes D_{<d}]_t.$ 

\begin{lemma}\label{ml2ax} Consider $r_i\in k_P^{\times}, \ 1\leq i\leq 5$ as in Lemma \ref{ml2} and $\lambda_i\in k(t), \ 1\leq i\leq 5$ as in Lemma \ref{Gauss}.
The following equality holds in $gr_d^{\mathcal{F}} \left [B_2(k(t))\otimes k(t)^{\times} \right ]_t:$
$$
(\{r_i\}_2 -\{\lambda_i\}_2)\otimes f_P=-\left\{\frac{r_i}{\lambda_i}  \right\}_2 \otimes (r_i-1)+\left\{\frac{r_i-1}{\lambda_i-1}  \right\}_2 \otimes \lambda_i.
$$
\end{lemma}

\begin{proof}
Thanks to the five-term relation for every $1\leq i \leq 5$
\be
\label{id19}
(\{r_i\}_2 -\{\lambda_i\}_2)\otimes f_P=-\left\{\frac{\lambda_i}{r_i}  \right\}_2 \otimes f_P-\left\{\frac{1-\lambda_i^{-1}}{1-r_i^{-1}}  \right\}_2 \otimes f_P-\left\{\frac{1-r_i}{1-\lambda_i}  \right\}_2 \otimes f_P.
\ee
We will simplify each of the three terms in the RHS of (\ref{id19}) in $gr_d^{\mathcal{F}}\left [ B_2(k(t)) \otimes k(t)^{\times}\right ]_t.$
Suppose that 
$$
h_i r_i =q_if_P+g_i.
$$ 
Since the degrees of $g_i$ and $h_i$ are less than $d,$ the degree of $q_i$ is also less than $d.$ Then,
\be
\begin{split}
	&-\left\{\frac{\lambda_i}{r_i}  \right\}_2 \otimes f_P=-\left\{\frac{g_i}{h_ir_i}  \right\}_2 \otimes f_P=\left\{\frac{h_ir_i}{g_i}  \right\}_2 \otimes f_P=\\
	&\left\{\frac{h_ir_i}{g_i}  \right\}_2 \otimes \left(\frac{h_ir_i}{g_i}-1 \right )-\left\{\frac{h_ir_i}{g_i}  \right\}_2 \otimes \frac{q_i}{g_i}=-\left\{\frac{r_i}{\lambda_i}  \right\}_2 \otimes \frac{q_i}{g_i}.\\
\end{split}
\ee
 We used here the following equalities in $[B_2(F)\otimes F^{\times}]_t$: $$\{x\}_2\otimes(x-1)=\{x\}_2\otimes(1-x)=-\{1-x\}_2\otimes(1-x)=0.$$ Similarly,
$$
\left\{\frac{1-\lambda_i^{-1}}{1-r_i^{-1}}  \right\}_2 \otimes f_P=-\left\{\frac{1-\lambda_i^{-1}}{1-r_i^{-1}}  \right\}_2 \otimes \frac{q_i}{(r_i-1)g_i}\\
$$
and
$$
	\left\{\frac{1-r_i}{1-\lambda_i}  \right\}_2 \otimes f_P=-\left\{\frac{1-r_i}{1-\lambda_i}  \right\}_2\otimes  \frac{q_i}{g_i-h_i}.
$$
So, we get the following:
$$
(\{r_i\}_2 -\{\lambda_i\}_2)\otimes f_P=-\left\{\frac{r_i}{\lambda_i}  \right\}_2 \otimes \frac{q_i}{g_i}+\left\{\frac{1-\lambda_i^{-1}}{1-r_i^{-1}}  \right\}_2 \otimes \frac{q_i}{(r_i-1)g_i}+\left\{\frac{1-r_i}{1-\lambda_i}  \right\}_2\otimes  \frac{q_i}{g_i-h_i}.
$$
From the five-term relation, we see that
\be
\begin{split}	
&-\left\{\frac{\lambda_i}{r_i}  \right\}_2 \otimes \frac{q_i}{(r_i-1)g_i}-\left\{\frac{1-\lambda_i^{-1}}{1-r_i^{-1}}  \right\}_2 \otimes \frac{q_i}{(r_i-1)g_i}-\left\{\frac{1-r_i}{1-\lambda_i}  \right\}_2 \otimes \frac{q_i}{(r_i-1)g_i}=\\
&(\{r_i\} -\{\lambda_i\})\otimes \frac{q_i}{(r_i-1)g_i}=0.\\
\end{split}
\ee
After adding it to the previous equality, we get that
$$
(\{r_i\}_2 -\{\lambda_i\}_2)\otimes f_P=-\left\{\frac{r_i}{\lambda_i}  \right\}_2 \otimes (r_i-1)+\left\{\frac{r_i-1}{\lambda_i-1}  \right\}_2 \otimes \frac{(r_i-1)\lambda_i}{\lambda_i-1},
$$
which finishes the proof of the statement of Lemma \ref{ml2ax}.
\end{proof}

\begin{proof}[Proof (End of Lemma \ref{ml2})]
Observe that the following relation holds in $B_P\subset gr_d^{\mathcal{F}} B_2(k(t)):$
\be\la{M1}
\left\{\frac{r_i-1}{\lambda_i-1}  \right\}_2+\left\{\frac{r_{i-1}r_{i+1}}{1-r_i}  \right\}_2-\left\{\frac{r_{i-1}}{\lambda_{i-1}}  \right\}_2-\left\{\frac{r_{i+1}}{\lambda_{i+1}} \right\}_2=0.
\ee
Indeed, the map $\delta$ applied to this expression lands in $\mathcal{F}_{d}\Lambda^2 k(t)^{\times}$
and vanishes there modulo terms of lower degree.

Our goal is to prove that 
$$
\sum_{i=1}^5 \Bigl (\{r_i\}_2\otimes f_P +\left \{ \frac{r_{i-1}r_{i+1}}{1-r_i}\right \}_2 \otimes r_i \Bigr )
$$
vanishes in $gr_d^{\mathcal{F}} \left [B_2(k(t))\otimes k(t)^{\times} \right ]_t.$ Using the $5-$term relation for $\lambda_i$ and the lemma above, we get that
\be\la{M2}
\begin{split}
	&\sum_{i=1}^5 \Bigl (\{r_i\}_2\otimes f_P +\left \{ \frac{r_{i-1}r_{i+1}}{1-r_i}\right \}_2 \otimes r_i \Bigr )=\\
	&\sum_{i=1}^5 \Bigl (\{r_i\}_2\otimes f_P -\{\lambda_i\}_2\otimes f_P  +\left \{ \frac{r_{i-1}r_{i+1}}{1-r_i}\right \}_2 \otimes r_i \Bigr )=\\
	&\sum_{i=1}^5 \Bigl (-\left\{\frac{r_i}{\lambda_i}  \right\}_2 \otimes (r_i-1)+\left\{\frac{r_i-1}{\lambda_i-1}  \right\}_2 \otimes \lambda_i  +\left \{ \frac{r_{i-1}r_{i+1}}{1-r_i}\right \}_2 \otimes \frac{r_i}{\lambda_i} +\left \{ \frac{r_{i-1}r_{i+1}}{1-r_i}\right \}_2 \otimes \lambda_i \Bigr ).\\
\end{split}
\ee
By Lemma \ref{ml1},
$$
\left \{ \frac{r_{i-1}r_{i+1}}{1-r_i}\right \}_2 \otimes \frac{r_i}{\lambda_i}+\left \{ \frac{r_i}{\lambda_i}\right \}_2 \otimes \frac{r_{i-1}r_{i+1}}{1-r_i}=0.
$$
We obtain 
\be
\begin{split}
	&\sum_{i=1}^5 \Bigl (\{r_i\}_2\otimes f_P +\left \{ \frac{r_{i-1}r_{i+1}}{1-r_i}\right \}_2 \otimes r_i \Bigr )=\\
	&\sum_{i=1}^5 \Bigl (\left\{\frac{r_i-1}{\lambda_i-1}  \right\}_2 \otimes \lambda_i  -\left \{ \frac{r_i}{\lambda_i}\right \}_2 \otimes (r_{i-1}r_{i+1})+\left \{ \frac{r_{i-1}r_{i+1}}{1-r_i}\right \}_2 \otimes \lambda_i \Bigr )=\\
	&\sum_{i=1}^5 \Bigl (\left\{\frac{r_{i-1}}{\lambda_{i-1}}  \right\}_2 \otimes \lambda_i+\left\{\frac{r_{i+1}}{\lambda_{i+1}} \right\}_2 \otimes \lambda_i-\left \{ \frac{r_i}{\lambda_i}\right \}_2 \otimes (r_{i-1}r_{i+1}) \Bigr )=\\
	&\sum_{i=1}^5 \Bigl (\left\{\frac{r_{i-1}}{\lambda_{i-1}}  \right\}_2 \otimes \frac{\lambda_i}{r_i} \Bigr ) + \sum_{i=1}^5 \Bigl (\left\{\frac{r_{i+1}}{\lambda_{i+1}} \right\}_2 \otimes \frac{\lambda_i}{r_i} \Bigr ).\\ 
\end{split}
\ee
The last expression vanishes, because by Lemma \ref{ml1}
$$
\left\{\frac{r_{i-1}}{\lambda_{i-1}}  \right\}_2 \otimes \frac{\lambda_i}{r_i} +\left\{\frac{r_{i}}{\lambda_{i}}  \right\}_2 \otimes \frac{\lambda_{i-1}}{r_{i-1}} =0.
$$
\end{proof}

\section{Proof of Lemma \ref{L1}}\label{S8}
For this we need to show that the chain map 
\be 
\begin{gathered}
    \xymatrix{
    & 0  \ar[r]& gr_1^{\mathcal{G}} gr_d^{\mathcal{F}}[ B_2(k(t)) \otimes \Lambda^{n-2} k(t)^{\times}]_t \ar[d] \ar[r] & gr_1^{\mathcal{G}} gr_d^{\mathcal{F}}[\Lambda^{n} k(t)^{\times} ]\ar[d]^{}\ar[r] & 0\\
    & 0 \ar[r]& \bigoplus_{\deg(P)=d} [B_2(k_P) \otimes \Lambda^{n-3} k_P^{\times}]_t \ar[r]& \bigoplus_{\deg(P)=d} \Lambda^{n-1} k_P^{\times} \ar[r]  &0}
 \end{gathered}
 \ee 
 is a quasi-isomorphism. That the chain map is a quasi-isomorphism when restricted to the second cohomology group, this follows from Theorem \ref{BT}. Denote by $K$ the kernel of the map
 $$
 gr_1^{\mathcal{G}} gr_d^{\mathcal{F}}[ B_2(k(t)) \otimes \Lambda^{n-2} k(t)^{\times}]_t \stackrel{\delta}{\lra} gr_1^{\mathcal{G}} gr_d^{\mathcal{F}}[\Lambda^{n} k(t)^{\times} ].
 $$
We must show that the induced map on the first cohomology group which is induced by the chain map $\oplus \partial_P$ is an isomorphism:
$$
\oplus_P \partial_P \colon K \lra \oplus_{\deg(P)=d} H^{n-2}_\mathrm{G}(k_P,\mathbb{Q}(n-1)).
$$

Recall that we defined the co-residue map $\sum c_P$   going in the opposite direction
from the chain map $\oplus \partial_P$. Since $\partial_P \circ c_P$ is the identity for $P=Q$ and vanishes otherwise, $\oplus \partial_P$ induces a surjective map on  cohomology. Let's prove that it is also injective. 
\begin{lemma}\label{L9}
	Vector space $gr_1^{\mathcal{G}} gr_d^{\mathcal{F}}[ B_2(k(t)) \otimes \Lambda^{n-2} k(t)^{\times}]_t$ is generated by elements 
	  \be
  \begin{split}
  	&1)\ \  \beta_{P}([r_1,r_2]_1)\otimes  g_{3} \wedge \ldots \wedge g_n,\\
  	&2)\ \  \left \{  \tilde{r} \right \}_2\otimes f_P  \wedge g_{4} \wedge \ldots \wedge g_n, \\
  \end{split}
  \ee
  where $r, r_1, r_2 \in k_P$ and $\deg(g_i)<d.$
\end{lemma}

\begin{proof}
 From  Corollary \ref{Cor} it follows that the space 
  $$
  gr_1^\mathcal{G}gr_d^{\mathcal{F}}\left [ B_2(k(t))\otimes \Lambda^{n-2} k(t)^{\times}\right ]_t
  $$
  is generated by elements of the following three types:
  \be
  \begin{split}
  	&a)\ \  \beta_{P}([r_1,r_2]_1)\otimes  g_{3} \wedge \ldots \wedge g_n,\\
  	&b)\ \  \left \{ \frac{g_1}{g_2} \right \}_2\otimes f_P \wedge g_{4} \wedge \ldots \wedge g_n, \\
  	&c)\ \  \beta_{g_1}([r_1,r_2]_1)\otimes f_P \wedge  g_{4} \wedge \ldots \wedge g_n, 
  \end{split}
  \ee
  where $g_1$ and $g_2$ are both irreducible polynomials of the same degree less than $d.$ 
  It is sufficient to consider the case $n=3.$ Our goal is to express elements of type $b)$ and $c)$ via elements of type $1)$ and $2).$ We will show that this is true for elements of type $b)$, type $c)$ is similar. Let $r$ be the residue of $\frac{g_1}{g_2}$ modulo $f_P.$ By the $5-$term relation
  $$
  \left \{ \frac{g_1}{g_2} \right \}_2 \otimes f_P-\left \{ r \right \}_2 \otimes f_P=\left (- \left \{ \frac{r  g_2}{g_1} \right \}_2 +\left \{ \frac{1-\frac{g_2}{g_1}}{1-\frac{1}{r}} \right \}_2 -\left \{ \frac{1-\frac{g_1}{g_2}}{1-r} \right \}_2 \right ) \otimes f_P
  $$ 
  the expression 
  $$
- \left \{ \frac{r  g_2}{g_1} \right \}_2 +\left \{ \frac{1-\frac{g_2}{g_1}}{1-\frac{1}{r}} \right \}_2 -\left \{ \frac{1-\frac{g_1}{g_2}}{1-r} \right \}_2 
  $$
  lies in $B_P,$ because $\delta$ sends it to $D_{<d}\otimes f_P$. 
  
  Finally, notice that for every $r_1, r_2 \in k_P,$ such that 
  $r_1 r_2=q  f_P+\overline{r_1r_2},$   
  $$
  \beta_{f_P}([r_1,r_2]_1)\otimes  f_P=\left \{ \frac{r_1  r_2}{\overline{r_1r_2}} \right \}_2\otimes  f_P =-\left \{ \frac{r_1  r_2}{\overline{r_1r_2}} \right \}_2\otimes  \frac{q}{\overline{r_1r_2}},$$
  because $\left \{ \dfrac{r_1  r_2}{\overline{r_1r_2}} \right \}_2\otimes  \left ( 1- \dfrac{r_1  r_2}{\overline{r_1r_2}} \right )=0.$ So, 
  $$\left \{ \frac{g_1}{g_2} \right \}_2 \otimes f_P-\left \{ r \right \}_2 \otimes f_P$$ is a linear combination of elements of type $1).$
\end{proof}

\begin{proof}[End of proof of Lemma \ref{L1}]
Let $x \in  gr_1^\mathcal{G}gr_d^{\mathcal{F}}\left [ B_2(k(t))\otimes \Lambda^{n-2} k(t)^{\times}\right ]_t$ be an element in the kernel of $\delta$ with vanishing residues $\partial_P$ at all points of degree $d.$ Then $x$ is 
lying in the subspace generated by elements of type $1)$ in the notation of Lemma \ref{L9}. Indeed, the map $id-\sum c_P \partial_P$ is well-defined and vanishes on elements of type 2). On the other hand,  $(id-\sum c_P \partial_P)x=x.$
It remains to notice that map $\delta$ is injective on the space of elements of type $1),$ which follows from the exact sequence
$$
 \bigoplus_{\deg(P)=d} \mathbb{S}^2 B_P \otimes \Lambda^{d-3} D_{<d} \lra  \bigoplus_{\deg(P)=d} B_P \otimes \Lambda^{d-2} D_{<d} \lra \bigoplus_{\deg(P)=d} f_P \otimes \Lambda^{d-1} D_{<d} \lra  \Lambda^{d-1} F_P^{\times}\lra 0
$$
and Lemma \ref{ml1}. 
\end{proof}

\section{Proof of Lemma \ref{L2}}\label{S9}
Here we prove Lemma \ref{L2}, claiming that   for $s>1$ the complex   
	$$
	gr_s^\mathcal{G}gr_d^{\mathcal{F}}{\cal B}_2(k(t),n)
	$$
	is exact. Recall that the space 
	$
	gr_s^\mathcal{G}gr_d^{\mathcal{F}} \Lambda^n k(t)^{\times}
	$
coincides with 
$\Lambda^s D_d \otimes   \Lambda^{n-s} D_{<d}.$
   Let $h$ be a map in the opposite direction
   $$
   h \colon gr_s^\mathcal{G}gr_d^{\mathcal{F}} \Lambda^n k(t)^{\times}
 \lra gr_s^\mathcal{G}gr_d^{\mathcal{F}}\left [ B_2(k(t))\otimes \Lambda^{n-2} k(t)^{\times}\right ]_t,
   $$
   sending 
   $$
   f_1 \wedge f_2 \wedge f_3 \wedge \ldots \wedge f_s \wedge g_{s+1}\wedge \ldots \wedge g_n
   $$  
   to
   $$
   -\left \{\frac{f_1}{f_2} \right \}\otimes f_3 \wedge \ldots \wedge f_s \wedge g_{s+1}\wedge \ldots \wedge g_n.
   $$
   Here we suppose that  $f_i$ are monic, irreducible, distinct polynomials of degree $d$ and the degrees of $g_i$ are less than $d.$
   
   \begin{lemma}
   Map $h$ is well-defined.	
   \end{lemma}
   \begin{proof}
   	It is sufficient to check that 
   	$$
   	\left \{\frac{f_1}{f_2} \right \}_2\otimes f_3 +\left \{\frac{f_3}{f_2} \right \}_2\otimes f_1
   	$$
   	 vanishes in $gr_3^\mathcal{G}gr_d^{\mathcal{F}}\left ( B_2(k(t))\otimes_a  k(t)^{\times}\right ).$
   	 From the five-term relation it follows that 
   	 $$
   	 \left \{\frac{f_1}{f_2} \right \}_2- \left \{\frac{f_3}{f_2} \right \}_2+\left \{\frac{f_3}{f_1} \right \}_2-
   	 \left \{\frac{1-\frac{f_2}{f_1}}{1-\frac{f_2}{f_3}} \right \}_2+\left \{\frac{1-\frac{f_1}{f_2}}{1-\frac{f_3}{f_2}} \right \}_2=0.
   	 $$
   	 After multiplication by $\dfrac{1-\frac{f_2}{f_1}}{1-\frac{f_2}{f_3}}=\dfrac{(f_1-f_2)f_3}{(f_3-f_2)f_1}$ we get the following equality in $B_2(k(t))\otimes_a k(t)^{\times}$:
   	 \be
   	 \begin{split}
   &	\left \{\frac{f_1}{f_2} \right \}_2\otimes f_3 +\left \{\frac{f_3}{f_2} \right \}_2\otimes f_1=\\
   &-\left \{\frac{f_1}{f_2} \right \}_2\otimes \frac{f_1-f_2}{(f_3-f_2)f_1}
	+ \left \{\frac{f_3}{f_2} \right \}_2\otimes \frac{(f_1-f_2)f_3}{f_3-f_2}-\\
	&\left \{\frac{f_3}{f_1} \right \}_2\otimes \frac{(f_1-f_2)f_3}{(f_3-f_2)f_1}-
   	 \left \{\frac{1-\frac{f_1}{f_2}}{1-\frac{f_3}{f_2}} \right \}_2\otimes \frac{(f_1-f_2)f_3}{(f_3-f_2)f_1}.
   	 \end{split}
   	 \ee
   	 All four expressions on the right side of this equality lie in $\mathcal{G}_{2}gr_d^{\mathcal{F}}\left [ B_2(k(t))\otimes \Lambda^{n-2} k(t)^{\times}\right ]_t.$
   \end{proof}
Obviously, $\delta \circ h=id,$ so to finish Lemma \ref{L2} it is sufficient to show that $h$ is surjective.
  From  Corollary \ref{Cor} it follows that the space 
  $$
  gr_s^\mathcal{G}gr_d^{\mathcal{F}}\left [ B_2(k(t))\otimes \Lambda^{n-2} k(t)^{\times}\right ]_t
  $$
  is generated by elements of the following three types:
  \be
  \begin{split}
  	&1)\ \  \beta_{f_1}([r_1,r_2]_1)\otimes f_2 \wedge\ldots \wedge f_s \wedge g_{s+1} \wedge \ldots \wedge g_n,\\
  	&2)\ \  \left \{ \frac{f_1}{f_2} \right \}_2\otimes f_3 \wedge\ldots \wedge f_s \wedge g_{s+1} \wedge \ldots \wedge g_n, \\
  	&3)\ \  B \otimes f_1 \wedge f_2 \wedge \ldots \wedge f_s \wedge g_{s+1} \wedge \ldots \wedge g_n, 
  \end{split}
  \ee
  where $f_1, \ldots f_s$ are monic, distinct, irreducible polynomials of degree $d,$ $g_{s+1}, \ldots, g_n$ are polynomials of degrees less than $d$ and $B$ lies in $\mathcal{F}_{<d}B_2(k(t)).$
  
  \begin{lemma}
  Elements of type 1) - 3) lie in the image of $h.$	
  \end{lemma}
  \begin{proof}
  For elements of type 2) there is nothing to prove.
  
 Let's prove it for elements of type 1). From Lemma \ref{L3} it follows that for any elements $X,Y \in B_2(F)$ we have
 $$
 X\otimes \delta (Y)-Y \otimes \delta (X) 
 $$
 vanishes in $[B_2(F) \otimes \Lambda^2 F^{\times}]_t.$ This can be applied to $X=B \in \mathcal{F}_{<d}B_2(k(x))$ and $Y=\left \{ \dfrac{f_1}{f_2} \right\}_2.$
 We get that 
 $$
 B \otimes \left (1-\frac{f_1}{f_2}\right ) \wedge \frac{f_1}{f_2} =\left \{ \frac{f_1}{f_2} \right \}_2 \otimes \delta(B).
 $$
 Since in 
  $
  gr_s^\mathcal{G}gr_d^{\mathcal{F}}\left [ B_2(k(t))\otimes \Lambda^{n-2} k(t)^{\times}\right ]_t
  $
  element $B \otimes f_1 \wedge f_2 \wedge \ldots \wedge f_s \wedge g_{s+1} \wedge \ldots \wedge g_n$ equals to $B \otimes \left (1-\frac{f_1}{f_2}\right ) \wedge \frac{f_1}{f_2} \wedge\ldots \wedge f_s \wedge g_{s+1} \wedge \ldots \wedge g_n,$ elements of type 1) lie in the image of $h.$
  
  It remains to deal with elements of type 3). It is enough to consider the case $s=2.$ Let $r$ be the residue of $\dfrac{\widetilde{r_1}\widetilde{r_2}}{\widetilde{r_1r_2}}$ modulo $f_2$. From the five term relation we see that 
  \be
  \begin{split}
  &\left \{ \dfrac{\widetilde{r_1}\widetilde{r_2}}{\widetilde{r_1r_2}} \right \}_2\otimes f_2=\left \{ r \right \}_2\otimes f_2-\left \{ \dfrac{\widetilde{r_1r_2} r}{\widetilde{r_1}\widetilde{r_2}} \right \}_2\otimes f_2+\\
  &\left \{ \dfrac{(\widetilde{r_1}\widetilde{r_2}-\widetilde{r_1r_2}) r}{\widetilde{r_1}\widetilde{r_2}(r-1)} \right \}_2\otimes f_2-
  \left \{ \dfrac{\widetilde{r_1}\widetilde{r_2}-\widetilde{r_1r_2}}{\widetilde{r_1r_2}(r-1)} \right \}_2\otimes f_2.\\
  \end{split}
  \ee
  All terms in the right side of the equality lie in the image of $h.$ We will show it for the second term, for the other three the proof is similar.
  Since $\deg(r)<d,$ there exist a polynomial $q$ of degree less than $d,$ such that 
  $$
  \widetilde{r_1}\widetilde{r_2}-r \widetilde{r_1r_2}=q f_2.
  $$
  So in $B_2(k(t))\otimes_a k(t)^{\times}$ we have
  \be
  \begin{split}
  &\left \{ \dfrac{(\widetilde{r_1}\widetilde{r_2}-\widetilde{r_1r_2}) r}{\widetilde{r_1}\widetilde{r_2}(r-1)} \right \}_2\otimes f_2=
  -\left \{ 1-\dfrac{(\widetilde{r_1}\widetilde{r_2}-\widetilde{r_1r_2}) r}{\widetilde{r_1}\widetilde{r_2}(r-1)} \right \}_2\otimes f_2\\
  &-\left \{ \dfrac{\widetilde{r_1r_2} r-\widetilde{r_1}\widetilde{r_2}}{\widetilde{r_1}\widetilde{r_2}(r-1)} \right \}_2\otimes f_2=
  -\left \{ \dfrac{-qf_2}{\widetilde{r_1}\widetilde{r_2}(r-1)} \right \}_2\otimes f_2=\left \{ \dfrac{-qf_2}{\widetilde{r_1}\widetilde{r_2}(r-1)} \right \}_2\otimes \dfrac{-q}{\widetilde{r_1}\widetilde{r_2}(r-1)}.
 \end{split}
  \ee
  Clearly $\delta\left \{ \dfrac{-qf_2}{\widetilde{r_1}\widetilde{r_2}(r-1)} \right \}_2-f_1 \wedge f_2 \in D_d \otimes D_{<d},$
  so by Corollary \ref{Cor}
  $$
  \left \{ \dfrac{-qf_2}{\widetilde{r_1}\widetilde{r_2}(r-1)} \right \}_2-\left \{\frac{f_1}{f_2}\right \}_2 
  $$
  lies in $\mathcal{G}_1 gr_d^\mathcal{F}B_2(k(t)).$
  So
  $$
  \left \{ \dfrac{(\widetilde{r_1}\widetilde{r_2}-\widetilde{r_1r_2}) r}{\widetilde{r_1}\widetilde{r_2}(r-1)} \right \}_2\otimes f_2-\left \{\frac{f_1}{f_2}\right \}_2 \otimes  \dfrac{-q}{\widetilde{r_1}\widetilde{r_2}(r-1)}
  $$
  vanishes in $gr_2^\mathcal{G}gr_d^\mathcal{F} \left [B_2(k(t))\otimes k(t)^{\times} \right]_t.$ 
  \end{proof}

\end{document}